\documentclass[11pt]{article}

\usepackage{latexsym,amsmath,amsthm,amsfonts,amsbsy,amscd,graphics,epsf}
\newtheorem{lemma}{Lemma}
\newtheorem{corollary}{Corollary}
\newtheorem{theorem}{Theorem}
\newtheorem{definition}{Definition}

\numberwithin{equation}{section}
\numberwithin{lemma}{section}
\numberwithin{corollary}{section}
\numberwithin{definition}{section}
\numberwithin{theorem}{section}

\newcommand{\qbin}[2]{\genfrac{[}{]}{0pt}{}{#1}{#2}}
\newcommand{\qbint}[2]{\left[\genfrac{}{}{0pt}{1}{#1}{#2}\right]}
\renewcommand{\pmod}[1]{\; (\bmod\, #1)}
\def\Integer{\mathbb{Z}}
\def\Rational{\mathbb{R}}
\def\bs{\boldsymbol}
\def\vL{\bs L}
\def\vl{\bs \ell}
\def\vm{\bs m}
\def\ve{{\bs e}}
\def\vx{\bs x}
\def\l{\ell}
\def\vn{\bs n}
\def\vA{\bs A}
\def\vu{\bs u}
\def\vQ{\bs Q}
\def\I{{\cal I}}
\def\vM{\bs M}
\def\by{\bar{y}}
\def\bl{\bar{l}}

\begin{document}

\title{Supernomial coefficients, polynomial identities and $q$-series}

\author{\Large
Anne Schilling
\thanks{Present address: Instituut voor Theoretische Fysica,
Universiteit van Amsterdam, Valckenierstraat 65, 1018 XE Amsterdam, 
The Netherlands}
\thanks{e-mail: {\tt schillin@phys.uva.nl}}\\
\mbox{} \\
\em Institute for Theoretical Physics, State University of New York \\
\em Stony Brook, NY 11794-3840, USA \\
\mbox{} \\
\Large and \\
\mbox{} \\
\setcounter{footnote}{0}
\Large S.~Ole Warnaar\footnotemark\setcounter{footnote}{2}
\thanks{e-mail: {\tt warnaar@phys.uva.nl}}\\
\mbox{} \\
\em Department of Mathematics, University of Melbourne\\
\em Parkville, Victoria 3052, Australia}
 
\date{\Large January, 1997}
\maketitle

\begin{abstract}
$q$-Analogues of the coefficients of $x^a$ in the expansion of
$\prod_{j=1}^N (1+x+\cdots +x^j)^{L_j}$ are proposed.
Useful properties, such as recursion relations, symmetries and
limiting theorems of the ``$q$-supernomial coefficients''
are derived, and a combinatorial interpretation using
generalized Durfee dissection partitions is given.
Polynomial identities of boson--fermion-type, based on
the continued fraction expansion of $p/k$ and involving the
$q$-supernomial coefficients, are proven. These include polynomial 
analogues of the Andrews--Gordon identities. 
Our identities unify and extend many of the 
known boson--fermion identities for 
one-dimensional configuration sums of solvable lattice models,
by introducing multiple finitization parameters.
\end{abstract}
 
\section{Introduction}
The Gaussian polynomial or $q$-binomial coefficient is defined as
\begin{equation}\label{qbin}
\qbin{L}{a}= \begin{cases} \cfrac{(q^{L-a+1})_{a}}{(q)_a} &
\text{for $a\in\Integer_+, L\in\Integer$},\\[3mm] 
0 & \text{otherwise},  \end{cases}
\end{equation}
where $\Integer_+$ is the set of non-negative integers,
$(x)_{\infty}=\prod_{i=0}^{\infty}(1-xq^i)$ and
$(x)_n=(x)_{\infty}/(xq^n)_{\infty}$ for $n\in\Integer$
(see for example~\cite{GR90}). 
Many nice identities involve the $q$-binomial, one of the simplest being
\begin{equation}\label{qbin_eqn}
(-x)_L = \sum_{a=0}^{\infty} \qbin{L}{a} x^a q^{a(a-1)/2},
\end{equation}
which is the $q$-analogue of the binomial expansion for $(1+x)^L$.
For non-negative $L$ equation~\eqref{qbin_eqn} asserts the combinatorial 
statement that
$q^{a(a-1)/2} \qbin{L}{a}$ is the generating function of
partitions with at most $a$ parts, no part exceeding $L-1$,
and all parts distinct. This result naturally leads to the question 
whether the multinomial coefficients defined by
\begin{equation}\label{trin}
(1+x+\cdots+x^N)^L =
\sum_{a=-\frac{NL}{2}}^{\infty}
\binom{L}{a}_{\! N} x^{a+\frac{NL}{2}},
\end{equation}
also admit a $q$-analogue.
(Note our convention that $a$ is half an odd integer when $NL$ is odd.)

For $N=2$ this question was addressed by Andrews and 
Baxter~\cite{AB87} who defined $q$-trinomial coefficients and studied 
some of their properties.  Interestingly though, equation~\eqref{trin} 
does not admit a $q$-analogue, and the combinatorial significance of 
the $q$-trinomial coefficients was only subsequently revealed, leading 
to a new proof of Schur's partition theorem~\cite{Andrews94}.
 
Recently, $q$-multinomial coefficients for arbitrary $N$ were 
introduced~\cite{S96,W96}, generalizing the Andrews--Baxter expressions. 
Surprisingly, the combinatorial interpretation of these $q$-multinomials
generalizes neither of the just mentioned partition results, but the fact 
that $q^{a^2} \qbin{L}{a}$ is the generating function of partions with 
exactly $a$ parts, no part exceeding $L$ and no part less than $a$.
In refs.~\cite{S96,W96} polynomial identities involving the $q$-multinomials 
were proven.  These identities can be viewed as finite or polynomial 
analogues of Rogers--Ramanujan or boson--fermion identities (see 
section~\ref{sec52} for the definition of the latter).  That is, for each 
value of a ``finitization parameter'' $L$ there is a polynomial identity, 
such that in the limit $L\to\infty$ a $q$-series identity of 
Rogers--Ramanujan or boson--fermion-type is recovered.

In this paper we generalize the $q$-multinomial coefficients
to ``$q$-supernomial coefficients'', or in short $q$-supernomials, 
by defining the natural 
$q$-analogues of the coefficients of $x^a$ in the expansion of
$\prod_{j=1}^N (1+x+\cdots+x^j)^{L_j}$. 
Using the $q$-supernomials, we prove polynomial identities labelled 
by not one but multiple finitization parameters.
First, a family of polynomial identities related to
Andrews' analytic counterpart of the Gordon partition theorem is considered. 
Subsequently, a much more general class of supernomial identities is proven, 
based on the continued fraction expansion of $p/k$. 
For both cases the proof relies on simple $q$-supernomial
recurrences and on known polynomial identities with a single finitization 
parameter, acting as initial conditions.  For the Andrews--Gordon case 
these are the polynomial boson--fermion identities of Foda and 
Quano~\cite{FQ95} and Kirillov~\cite{Kirillov95}. 
For the continued fraction identities, the initial conditions are the 
polynomial identities of Berkovich and McCoy~\cite{BM96,BMS96} for 
finite analogues of the minimal Virasoro characters $\chi^{(p-k,p)}$.
Besides the special cases used as initial conditions,
our supernomial identities take many of the known 
boson--fermion polynomial identities for one-dimensional configuration 
sums of solvable lattice models as special one-variable subcases.
Taking the $q$-series limits of the $q$-supernomial identities,
we obtain new boson--fermion type identities for generalizations of the
A$^{(1)}_1$ branching functions.

The first part of this paper, comprised of sections~\ref{sec2} and
\ref{sec3}, develops the general theory of $q$-supernomials.
Section~\ref{sec2} is devoted to the definition of the $q$-supernomials as 
$N$-variable generalizations of the $q$-multinomial coefficients,
and to stating its most fundamental properties, such as symmetries, 
recurrences and $q$-series limits. In section~\ref{sec3} a combinatorial 
interpretation of $q$-supernomials is given using Durfee dissections. 
This interpretation is then employed to show how several of the results 
established in section~\ref{sec2} can be understood purely combinatorially.
The second part of this paper deals with the application of
$q$-supernomials to polynomial identities labelled by multiple finitization 
parameters. The supernomial identities related to the Andrews--Gordon 
identity are given in section~\ref{sec4}. The supernomial identities based
on the continued fraction expansion of $p/k$, and their corresponding 
$q$-series identities, are the subject of section~\ref{sec5}. 
We conclude with an outlook on open problems related
to the $q$-supernomials.

\section{$q$-Supernomial coefficients}\label{sec2}
\subsection{Preliminaries}\label{sec_prelims}
Before defining the $q$-supernomial coefficients, we
introduce some general notation used throughout this paper.
$q$-Supernomials are multivariable generalizations of the 
$q$-multinomial coefficients.
It is therefore useful to introduce a vector notation, and
for a fixed positive integer $N$, we set $\vL = (L_1,\dots,L_N)$, where
each of the entries $L_j\in\Integer$.
We further need the Cartan matrix $C$ of the Lie algebra A$_{N-1}$, with
components $C_{i,j}=2\delta_{i,j}-\delta_{|i-j|,1}$, and
the Cartan-type matrix $T$, associated with the tadpole graph of $N$ nodes,
with components 
$T_{i,j}=\delta_{i,j} (2-\delta_{j,N})-\delta_{|i-j|,1}$.
The entries of the inverse of these matrices are given by
\begin{equation}\label{CT_inv}
\begin{array}{ll}
C_{i,j}^{-1} = \min\{i,j\}-\dfrac{ij}{N} \quad &  i,j=1,\dots,N-1, \\[2mm]
T_{i,j}^{-1} = \min\{i,j\} &  i,j=1,\dots,N. 
\end{array}
\end{equation}
Besides $\vL$, it is often convenient to use the vector $\vl = T^{-1} \vL$.
(Throughout this paper transposition symbols are omitted.)
Finally we use the unit vectors $\ve_n$ with entries 
$(\ve_n)_i = \delta_{n,i}$. The dimension of $\ve_n$ will be either explicitly
stated or will be clear from the context. 
For $N$-dimensional unit vectors we use the convention that 
$\ve_n={\bs 0}$ for $n\neq 1,\ldots,N$.

For the $q$-binomial \eqref{qbin} we need the transformation
\begin{equation}\label{dual_qbin}
\qbin{L}{a}_{1/q}=q^{-a(L-a)}\qbin{L}{a},
\end{equation}
as well as the recurrences
\begin{align}\label{bi1}
\qbin{L}{a}&= \qbin{L-1}{a-1}+q^a \qbin{L-1}{a}, \\[2mm]
\label{bi2}
\qbin{L}{a}
&=\qbin{L-1}{a}+q^{L-a} \qbin{L-1}{a-1}.
\end{align}

\subsection{Supernomial coefficients}
We are now prepared to introduce the supernomial coefficients. 
Before introducing a base, the $q=1$ case is briefly considered first.
The supernomial coefficients are defined as the following generalization 
of the multinomial coefficients
\begin{definition}\label{def_supernomials}
Let $N$ be a positive integer, $\vL\in \Integer^N$ and $\vl = T^{-1} \vL$.
For $a+\frac{1}{2}\ell_N\in\Integer_+$
we define the supernomial coefficient $\binom{\vL}{a}$ through
\begin{equation}\label{supdef}
\prod_{j=1}^N (1+x+\cdots +x^j)^{L_j} = 
\sum_{a=-\frac{\ell_N}{2}}^{\infty}
\binom{\vL}{a} x^{a+\frac{\ell_N}{2}}.
\end{equation}
\end{definition}
Clearly, for $\vL=L \ve_n$ the supernomial simplifies to the 
($n$-)multinomial coefficient~\eqref{trin}, 
$\binom{L\ve_n}{a}=\binom{L}{a}_n$.
Some simple properties of the supernomials are the initial 
condition and symmetry
\begin{equation}\label{supsymm}
\binom{\bs 0}{a} =  \delta_{a,0} \qquad  \text{and}\qquad
\binom{\vL}{a} =  \binom{\vL}{-a} \quad \text{for $\vL\in\Integer^N_+$},
\end{equation}
as well as the recurrences
\begin{equation}\label{recurrences}
\binom{\vL}{a} = \binom{\vL+\ve_{n-1}-2\ve_n+\ve_{n+1}}{a} +
\binom{\vL-2\ve_n}{a} \qquad \text{for $n=1,\dots,N-1$}.
\end{equation}
To see this last result we note the identity
\begin{equation*}
(1+x+\cdots+x^n)^2 = 
(1+x+\cdots+x^{n-1})(1+x+\cdots+x^{n+1}) + x^n.
\end{equation*}

A representation useful for defining the $q$-analogue of the
supernomials is provided by
\begin{equation}\label{supernomial}
\binom{\vL}{a} =
\sum_{j_1+\cdots+j_N=a+\frac{\ell_N}{2}}
\binom{L_N}{j_N}\binom{L_{N-1}+j_N}{j_{N-1}}\cdots\binom{L_1+j_2}{j_1}.
\end{equation}
This follows by $N$-fold use of the binomial expansion,
\begin{align}
\sum_{a=-\frac{\ell_N}{2}}^{\infty}
\binom{\vL}{a} x^{a+\frac{\ell_N}{2}} & =
\sum_{j_1,\dots,j_N \geq 0} x^{j_1+\dots+j_N}
\prod_{k=1}^N \binom{L_k+j_{k+1}}{j_k} \qquad \qquad (j_{N+1}=0) \notag \\
& = (1+x)^{L_1} \sum_{j_2,\dots,j_N} x^{j_2+\dots+j_N}(1+x)^{j_2}
\prod_{k=2}^N \binom{L_k+j_{k+1}}{j_k} \notag \\
& = (1+x)^{L_1} (1+x+x^2)^{L_2}
\sum_{j_3,\dots,j_N} x^{j_3+\dots+j_N}(1+x+x^2)^{j_3}
\prod_{k=3}^N \binom{L_k+j_{k+1}}{j_k} \notag \\
& = \dots = 
\prod_{j=1}^N (1+x+\cdots +x^j)^{L_j}. \notag
\end{align}

\subsection{$q$-Supernomial coefficients}
We propose the following $q$-analogue of~\eqref{supernomial}.
\begin{definition}\label{sdef1}
Let $N$ be a positive integer, $\vL\in \Integer^N$ and $\vl = T^{-1} \vL$.
Then for $a+\frac{1}{2}\ell_N\in\Integer_+$
\begin{equation}\label{qsupernomial1}
\qbin{\vL}{a} =\sum_{j_1+\cdots+j_N=a+\frac{\ell_N}{2}}
q^{\sum_{k=2}^N j_{k-1} (L_k+\cdots+L_N-j_k)}
\qbin{L_N}{j_N}\qbin{L_{N-1}+j_N}{j_{N-1}}\cdots\qbin{L_1+j_2}{j_1}.
\end{equation}
\end{definition}

We also need the $q\to 1/q$ form of~\eqref{qsupernomial1}.
\begin{definition}
With the same parameters as in definition~\ref{sdef1} the
$q$-supernomial $T(\vL,a)$ is defined as
\begin{equation}\label{duality}
T(\vL,a)=q^{\frac{1}{4}\vL T^{-1}\vL-\frac{a^2}{N}} \qbin{\vL}{a}_{1/q}.
\end{equation}
\end{definition}

Since $\qbin{(L_1,\dots,L_N,0,\dots,0)}{a} = \qbin{(L_1,\dots,L_N)}{a}$
this gives
\begin{equation}\label{T_M}
T\bigl((L_1,\dots,L_N,0,\dots,0),a\bigr)
=q^{\frac{M-N}{MN}a^2}T\bigl((L_1,\dots,L_N),a\bigr),
\end{equation}
where the vector $(L_1,\dots,L_N,0,\dots,0)$ is $M$-dimensional.
Using \eqref{dual_qbin} we also deduce that
\begin{equation}\label{T_L1}
T(L\ve_1,a)=q^{\frac{N-1}{N}a^2}\qbin{L}{\frac{L}{2}+a},
\end{equation}
where $\ve_1$ is $N$-dimensional.
In section~\ref{sec5} these equations will play quite a remarkable role.

Let us now state an explicit formula for $T(\vL,a)$ when $\vL\in\Integer_+^N$.
\begin{lemma}\label{lem_explicit}
Let $N$ be a positive integer, $\vL\in \Integer^N_+$ and $\vl = T^{-1} \vL$.
Set $L_0=0$,
\begin{equation}\label{xjdef}
x_j = q^{\sum_{k=0}^j \left(m_k-\frac{L_k}{2}\right)}\qquad (j=1,\ldots,N),
\end{equation}
and 
\begin{equation}\label{eqn_restr}
m_0=\frac{\ell_1}{2}+\frac{a}{N}-(C^{-1}\vm)_1, \qquad
m_N=-\frac{a}{N}-(C^{-1}\vm)_{N-1},
\end{equation}
where $\vm=(m_1,\dots,m_{N-1})$.
Then for $a=-\frac{1}{2}\ell_N,\dots,\frac{1}{2}\ell_N$,
\begin{equation}\label{qsupernomial2}
T(\vL,a)=\sideset{}{'}
\sum_{\substack{m_i \in \Integer_+ - \frac{1}{2}L_i \\[0.5mm]
0\leq i \leq N}}
q^{\vm C^{-1}\vm} \prod_{j=0}^N
\frac{(x_j q)_{L_j}}{(q)_{\frac{L_j}{2}+m_j}} \: ,
\end{equation}
where the primed summation symbol denotes a 
sum over $m_0,\dots,m_N$, such that \eqref{eqn_restr} is satisfied.
\end{lemma}
In actual use of $T$ it is helpful to note that $x_N$ defined in 
\eqref{xjdef} is equal to one. This can be simply observed by noting 
that $(C^{-1} \vm)_1 + (C^{-1} \vm)_{N-1} = m_1 + \cdots + m_{N-1}$.

\begin{proof}[Proof of lemma \ref{lem_explicit}]
In the following we set $j_{N+1}=0$.
Using \eqref{dual_qbin} one may derive 
\begin{equation*}
\qbin{\vL}{a}_{1/q}=
\sum_{j_1+\cdots+j_N=a+\frac{\ell_N}{2}}
q^{\sum_{k=1}^N j_k(j_k-L_{k}-\cdots-L_N)}
\prod_{k=1}^N \qbin{L_k+j_{k+1}}{j_k}.
\end{equation*}
Rewriting the $q$-binomials using $(q)_{m+n} = (q)_m(q^{1+m})_n$, we obtain
\begin{equation*}
\prod_{k=1}^N \qbin{L_k+j_{k+1}}{j_k}
=\frac{1}{(q)_{j_1}}
\prod_{k=1}^{N} \frac{(q^{1+j_{k+1}})_{L_k}}{(q)_{L_k+j_{k+1}-j_k}}
=\prod_{k=0}^{N} \frac{(x_k\,q)_{L_k}}
{(q)_{\frac{L_k}{2}+m_k}}.
\end{equation*}
Notice that we used $L_k\geq 0$ to derive the first equality since the 
rewriting $(q^{L+1})_j/(q)_j=(q^{j+1})_L/(q)_L$ is problematic when $L<0$.
The last step follows from $L_0=0$ and the variable change
\begin{equation}\label{eqn_jk}
j_k=\sum_{i=0}^{k-1}\left(m_i-\frac{L_i}{2}\right), \quad 
1\leq k\leq N+1,
\end{equation}
where $m_i\in\Integer-\frac{L_i}{2}$ $(i=0,\ldots,N)$, and $m_0$ and 
$m_N$ satisfy~\eqref{eqn_restr}.
This variable change makes sense because one can check that 
$j_{N+1}=0$, $j_1+\cdots+j_N=a+\frac{\ell_N}{2}$ and 
$j_k\in \Integer,~(1\leq k\leq N)$. 
The equality of the exponents
\begin{equation*}
\frac{1}{4}\vL T^{-1} \vL-\frac{a^2}{N}
+\sum_{k=1}^N j_k(j_k-L_{k}-\cdots-L_N)=\vm C^{-1}\vm
\end{equation*}
follows from (\ref{eqn_jk}) and the explicit form of the inverse
Cartan matrices given in (\ref{CT_inv}).
\end{proof}

We note that the $q$-supernomials for $\vL=L\ve_n$ $(L\geq 0)$ and the 
$q$-multinomials of refs.~\cite{S96,W96} with label $\ell=0$ coincide,
$$\qbin{L\ve_n}{a-\frac{nL}{2}}=\qbin{L}{a}^{(0)}_n \qquad \text{and}
\qquad T(L\ve_n,a)=q^{\frac{N-n}{Nn}a^2}T_0^{(n)}(L,a).$$
For general $\ell$, the $q$-multinomials are, however, only equal to the 
$q$-supernomials in the limit $L\to\infty$ 
(see equation~\eqref{multi_limit} below).

\subsection{Symmetries and recurrences of the $q$-supernomial coefficients}
In this section we prove the $q$-analogues of the supernomial 
symmetry~\eqref{supsymm} and recurrences~\eqref{recurrences}. The recurrences 
play a crucial role in the proof of the polynomial identities in 
sections~\ref{sec4} and~\ref{sec5}.

\begin{lemma}\label{lem_sym}
For $\vL\in\Integer^N_+$, the $q$-supernomials satisfy the symmetries
\begin{equation}\label{sym}
\qbin{\vL}{a}=\qbin{\vL}{-a} \qquad \text{ and } \qquad 
T(\vL,a)=T(\vL,-a).
\end{equation}
\end{lemma}

\begin{lemma}\label{lem_rec}
For $1\leq n\leq N-1$ and $\vL\in\Integer^N$, the following recursion 
relations hold
\begin{align}\label{eqn_rec}
\qbin{\vL}{a}&=q^{\l_n-n}
\qbin{\vL-2\ve_n}{a}+\qbin{\vL+\ve_{n-1}-2\ve_{n}+\ve_{n+1}}{a},\\[4mm]
\label{eqn_Trec}
T(\vL,a)&=T(\vL-2\ve_n,a)
+q^{\frac{1}{2}(L_n-1)}T(\vL+\ve_{n-1}-2\ve_n+\ve_{n+1},a).
\end{align}
\end{lemma}

The next lemma states which initial conditions are sufficient to determine
a function satisfying the $q$-supernomial recurrences.

\begin{lemma}\label{lem_complete}
Let $X$ be a function of $\vL\in \Integer^N$ which obeys the
recurrences
\begin{equation}\label{Xrec}
X(\vL)=q^{\l_n-n}X(\vL-2\ve_n)+X(\vL+\ve_{n-1}-2\ve_n+\ve_{n+1})
\end{equation}
for $1\leq n\leq N-1$ and all $\vL\in\Integer^N$. Then $X(\vL)$ for 
$\vL\in\Integer^N_+$ is uniquely determined by $X(L\ve_1)$ with $L\geq 0$.
\end{lemma}

We now give an analytic proof of  lemma~\ref{lem_rec}. In the next section,
where a partition theoretical interpretation of the $q$-supernomials will
be discussed, we give an alternative, combinatorial proof.

\begin{proof}[Proof of lemma \ref{lem_rec}]
Equation \eqref{eqn_Trec} follows immediately from \eqref{eqn_rec} 
using~\eqref{duality}. 

For the proof of \eqref{eqn_rec} we repeatedly use
the $q$-binomial recurrences \eqref{bi1} and \eqref{bi2}.
After inserting~\eqref{qsupernomial1} into the right-hand side 
of~\eqref{eqn_rec}, we change variables 
$j_i\to j_i-1$ for $1\leq i\leq n$ in the first term. In the second term we
apply \eqref{bi2} to the $q$-binomial containing $L_{n+1}$. This yields
for the right-hand side of \eqref{eqn_rec} (with $j_0=j_{N+1}=0$),
\begin{multline}
\sum_{j_1+\dots+j_N=a+\frac{\l_N}{2}} \!
q^{\sum_{k=2}^N j_{k-1} (L_k+\cdots+L_N-j_k)}  
\Bigl( q^{L_1+\cdots+L_n-j_1+j_n+j_{n+1}-1}
\prod_{k=1}^N \qbint{L_k+j_{k+1}-\delta_{k,n}-\theta(k\leq n)}
{j_k-\theta(k\leq n)} \\ 
+q^{j_n-j_{n-1}}
\prod_{k=1}^N \qbint{L_k+j_{k+1}+\delta_{k,n-1}-2\delta_{k,n}}{j_k}
+q^{L_{n+1}-j_{n-1}+j_n-j_{n+1}+j_{n+2}+1}\prod_{k=1}^N
\qbint{L_k+j_{k+1}+\delta_{k,n-1}-2\delta_{k,n}}{j_k-\delta_{k,n+1}}\Bigr),
\notag
\end{multline}
where $\theta(\text{true})=1$ and $\theta(\text{false})=0$.
Inserting the telescopic expansion
\begin{equation*}
\prod_{k=1}^{n-1} \qbint{L_k+j_{k+1}+\delta_{k,n-1}}{j_k}
=\sum_{m=0}^{n-1}q^{j_m}\prod_{k=1}^{n-1} \qbint{L_k+j_{k+1}+\delta_{k,n-1}
-\theta(m\leq k<n)}{j_k-\theta(m<k<n)}
\end{equation*}
in the second term, yields
\begin{multline}\label{eqn_rpc}
\sum_{j_1+\dots+j_N=a+\frac{\l_N}{2}} \!
q^{\sum_{k=2}^N j_{k-1} (L_k+\cdots+L_N-j_k)} 
\Bigl( q^{L_{n+1}-j_{n-1}+j_n-j_{n+1}+j_{n+2}+1}\prod_{k=1}^N
\qbint{L_k+j_{k+1}+\delta_{k,n-1}-2\delta_{k,n}}{j_k-\delta_{k,n+1}}\\
+\sum_{m=0}^{n-1} q^{j_m-j_{n-1}+j_n}\prod_{k=1}^N 
\qbint{L_k+j_{k+1}+\delta_{k,n-1}-2\delta_{k,n}-\theta(m\leq k<n)}
{j_k-\theta(m<k<n)}\\
+q^{L_1+\cdots+L_n-j_1+j_n+j_{n+1}-1}
\prod_{k=1}^N \qbint{L_k+j_{k+1}-\delta_{k,n}-\theta(k\leq n)}
{j_k-\theta(k\leq n)}\Bigr).
\end{multline}
Now change $j_{n+1}\to j_{n+1}+1$ and $j_n\to j_n-1$ in the first term
and $j_n\to j_n-1$ and $j_{m+1}\to j_{m+1}+1$ in the $m$-th term
in the sum over $m$. This leads to (changing $m\to m-1$ and including 
the last term in~\eqref{eqn_rpc} as the $m=0$ term in the sum) 
\begin{multline}
\sum_{j_1+\dots+j_N=a+\frac{\l_N}{2}} \!
q^{\sum_{k=2}^N j_{k-1} (L_k+\cdots+L_N-j_k)} \Bigl(
\prod_{k=1}^N \qbint{L_k+j_{k+1}-\delta_{k,n}}{j_k-\delta_{k,n}} \\
+ \sum_{m=0}^n q^{L_{m+1}+\cdots+L_n-j_{m+1}+j_n+j_{n+1}+\delta_{m,n}-1} 
\prod_{k=1}^N \qbint{L_k+j_{k+1}-\delta_{k,n}-\theta(m\leq k\leq n)}
{j_k-\theta(m<k\leq n)} \Bigr). \notag
\end{multline}
Now one can telescopically combine all the terms, starting with 
combining the $m=0$ and $m=1$ terms in the sum using (\ref{bi2}), 
then combining this with the $m=2$ term, and so on.
The result of this can be combined with the term in the first line 
using (\ref{bi1}). This yields \eqref{qsupernomial1} and we are done.
\end{proof}

\begin{proof}[Proof of lemma \ref{lem_complete}]
We will show that the recurrences
\begin{multline}\label{eqn_recs}
X(\vL+\ve_{n+1})=q^{-\l_{n-2}-2(n-2)}\Bigl(X(\vL+\ve_{n-2}+\ve_{n-1}
+\ve_n)-q^{\l_{n-1}+2n-3}X(\vL+\ve_{n-1})\\
-\theta(n>2) \bigl( q^{\l_n+2n-3}X(\vL+\ve_{n-3})
+X(\vL+\ve_{n-3}+\ve_{n-1}+\ve_{n+1}) \bigr) \Bigr)
\end{multline}
for $2\leq n\leq N-1$ and $L_{n-1}=0$ follow from the 
recurrences~\eqref{Xrec}.
These recurrences together with the recurrences~\eqref{Xrec}, both 
for $\vL\in\Integer^N_+$, uniquely determine
$X(\vL)$ for $\vL\in\Integer^N_+$ from the initial condition $X(L\ve_1)$
$(L\geq 0)$ for the following reason.
First, $X(L_1\ve_1+L_2\ve_2)$ for $L_1,L_2\geq 0$ is determined by
recurrence~\eqref{Xrec} for $n=1$. Now assume that $X(\vL')$ is known for
all $\vL'=L'_1\ve_1+\cdots+L'_{n+1}\ve_{n+1}$ $(n\geq 1)$ with 
$L'_i\in\Integer_+~(1\leq i\leq n)$ and $0\leq L'_{n+1}<L_{n+1}$ for some 
$L_{n+1}>1$. It then suffices to show that $X(\vL)$ with
$\vL=L_1\ve_1+\cdots+L_{n+1}\ve_{n+1}$ and $L_i\in\Integer_+~(1\leq i\leq n)$ 
is determined by the recurrences~\eqref{Xrec} and~\eqref{eqn_recs}.
Rewriting~\eqref{Xrec} as
\begin{equation}\label{rec_solved}
X(\vL+\ve_{n+1})=X(\vL-\ve_{n-1}+2\ve_n)-q^{\l_n+1}X(\vL-\ve_{n-1})
\end{equation}
we see that $X(\vL)$ with $L_{n-1}>0$ is determined from $X(\vL')$ with
only positive components $L_i'$. $X(\vL)$ with $L_{n-1}=0$ now follows
from~\eqref{eqn_recs} since all terms on the right-hand side are known
and have positive components.

It remains to be shown that~\eqref{eqn_recs} can indeed be deduced 
from~\eqref{Xrec} for $\vL\in\Integer^N$. We start with~\eqref{rec_solved}
and apply~\eqref{Xrec} with $n$ replaced by $n-1$ to the first term on the 
right-hand side
\begin{multline}
X(\vL+\ve_{n+1})=X(\vL-\ve_{n-2}+\ve_{n-1}+\ve_n)-q^{\l_{n-1}+1}
X(\vL-\ve_{n-2}-\ve_{n-1}+\ve_n)\\
-q^{\l_n+1}X(\vL-\ve_{n-1}).\notag
\end{multline}
Once more applying~\eqref{Xrec} with $n\to n-1$ now to the second term on the 
right-hand side yields
\begin{multline}\label{inter}
X(\vL+\ve_{n+1})=X(\vL-\ve_{n-2}+\ve_{n-1}+\ve_n)-q^{\l_{n-1}+1}
X(\vL-2\ve_{n-2}+\ve_{n-1})\\
+q^{2\l_{n-1}-2n+5}X(\vL-2\ve_{n-2}-\ve_{n-1})-q^{\l_n+1}X(\vL-\ve_{n-1}).
\end{multline}
Observe that for $n=2$ and $L_1=0$ the last two terms cancel yielding
\eqref{eqn_recs} for $n=2$. For $n>2$, the last two terms 
combine to $-q^{\l_n+1}X(\vL+\ve_{n-3}-2\ve_{n-2})$ by~\eqref{Xrec}
if $L_{n-1}=0$.
Replacing $X(\vL+\ve_{n+1})$ on the left-hand side of~\eqref{inter} by
\begin{equation*}
X(\vL+\ve_{n+1})=q^{\l_{n-2}}X(\vL-2\ve_{n-2}+\ve_{n+1})+
X(\vL+\ve_{n-3}-2\ve_{n-2}+\ve_{n-1}+\ve_{n+1}),
\end{equation*}
solving for $X(\vL-2\ve_{n-2}+\ve_{n+1})$ and replacing 
$\vL\to \vL+2\ve_{n-2}$ yields~\eqref{eqn_recs}.
\end{proof}

\begin{proof}[Proof of lemma \ref{lem_sym}]
Lemma~\ref{lem_sym} follows from lemmas~\ref{lem_rec} and~\ref{lem_complete}
by induction. Setting $X(\vL)=\qbint{\vL}{a}$ and recalling that
the $q$-supernomial fulfills the recurrence~\eqref{eqn_rec} we know from
lemma~\ref{lem_complete} and its proof that $\qbint{\vL}{a}$ for 
$\vL\in\Integer_+^N$ is determined from $\qbint{L\ve_1}{a}$ $(L\geq 0)$ from
recurrences only involving $\qbint{\vL}{a}$ with $\vL\in\Integer_+^N$.
Since $\qbin{L\ve_1}{a}=\qbin{L\ve_1}{-a}$ thanks to the symmetry of the
$q$-binomial the lemma is proven.
\end{proof}

\subsection{Limiting behaviour of the $q$-supernomials}
By definition both types of $q$-supernomials are polynomials in $q$. 
In the limit $L_m\to\infty$ the $q$-supernomials yield 
$q$-series, as described in the next few lemmas.
\begin{lemma}\label{lem_limit1}
For $m=1,2,\dots,N$ and $|q|<1$,
\begin{equation}\label{limit1}
\lim_{L_m\to\infty} \qbin{\vL}{a-\frac{\ell_N}{2}}=\frac{1}{(q)_a}\, ,
\end{equation}
independent of $m$.
\end{lemma}
\begin{proof}
First notice that for $|q|<1$ 
$$q^{\sum_{k=2}^m j_{k-1}\left( L_k+\cdots+L_N-j_k\right)}$$
is only non-zero for $L_m\to\infty$ if $j_1=\cdots=j_{m-1}=0$. This implies
\begin{equation*}
\lim_{L_m\to\infty} \qbin{\vL}{a-\frac{\ell_N}{2}}=
\sum_{j_m+\cdots+j_N=a} q^{\sum_{k=m+1}^N j_{k-1}\left( L_k+\cdots+L_N-j_k
\right)} \qbin{L_N}{j_N}
\cdots \qbin{L_{m+1}+j_{m+2}}{j_{m+1}}\frac{1}{(q)_{j_m}}.
\end{equation*}
Hence it is sufficient to show that 
\begin{equation*}
\lim_{L_1\to\infty} \qbin{\vL}{a-\frac{\ell_N}{2}}=\frac{1}{(q)_a}.
\end{equation*}
Equation \eqref{eqn_rec} implies that
\begin{equation}\label{lim_n}
\lim_{L_1\to\infty} \qbin{\vL}{a-\frac{\ell_N}{2}}=\lim_{L_1\to\infty}
\qbin{\vL+\ve_{n-1}-2\ve_n+\ve_{n+1}}{a-\frac{\ell_N}{2}}
\end{equation}
for $n=1,2,\dots,N-1$. From \eqref{lim_n} with $n=1$ it follows that
$\lim_{L_1\to\infty} \qbin{\vL}{a-\ell_N/2}$ 
is independent of $L_2$. From 
\eqref{lim_n} with $n=2$ one deduces that it is independent of $L_3$,
et cetera. We can therefore choose $L_2=\cdots=L_N=0$, which immediately
implies \eqref{limit1}.
\end{proof} 

In contrast to the result of lemma~\ref{lem_limit1},
the various limits of $T(\vL,a)$ depend on which
components of $\vL$ are taken to infinity, and to give our results we
need some definitions first.
\begin{definition}\label{def_b}
Let $N,h$ be integers such that $1\leq h\leq N$. Choose integers
$k_i$ $(1\leq i\leq h)$ such that $1\leq k_1<k_2<\cdots <k_h\leq N$ and
denote the sets $K=\{k_1,k_2,\dots,k_h\}$ and 
$\bar{K}=\{1,2,\dots,N\}-K$.
Finally, fix $a\in\Integer$, $L_k\in\Integer_+$ for $k\in\bar{K}$ and
$\sigma_{k}=0,1$ for $k\in K$ such that $a+\sum_{k\in\bar{K}} kL_k
+\sum_{k\in K}k\sigma_k$ is even. Then 
\begin{equation}\label{b_def}
b_a^{\{L_k|k\in\bar{K}\}\{\sigma_k|k\in K\}}(q)
=\sum_{\substack{m_i\in\Integer_+-\frac{1}{2}L_i,~i\in \bar{K}\\
m_i\in \Integer-\frac{1}{2}\sigma_i,~i\in K \\
\frac{a}{2N}+(C^{-1}\vm)_{N-1}=-m_N}}
q^{\vm C^{-1} \vm}\frac{(x_{k_h}q)_{\infty}}{(q)_{\infty}^{h+1}}
\frac{\prod_{j=k_h+1}^N (x_jq)_{L_j}}{\prod_{j\in \bar{K}}
(q)_{\frac{L_j}{2}+m_j}},
\end{equation}
with the variables $x_j$ given by
\begin{equation} \label{xj}
x_j=q^{\frac{a}{2N}+\frac{1}{2}(L_{j+1}+\cdots+L_N)+(\ve_j-\ve_{j+1})C^{-1}\vm}
\quad(k_h\leq j<N)\quad\text{and}\quad x_N=1.
\end{equation}
\end{definition}

\begin{definition}\label{def_super_string}
Let $N$ be a positive integer, $\vL\in\Integer_+^{N-1}$,
$a\in\Integer$ and $\sigma=0,1$
such that $r=a-N(C^{-1}\vL)_{N-1}+N\sigma$ is even. Then
\begin{equation}\label{super_string}
c_a^{\vL,\sigma}(q)=\frac{q^{\tfrac{\vL C^{-1}\vL}{2(N+2)}}}{(q)_{\infty}}
\sum_{\substack{\vm\in\Integer_+^{N-1}\\ 
\frac{r}{2N}-(C^{-1}\vm)_1\in\Integer}}
\frac{q^{\vm C^{-1}(\vm-\vL)}}{(q)_{m_1}\dots (q)_{m_{N-1}}}.
\end{equation}
\end{definition}

The function $c_a^{\vL,\sigma}(q)$ has the symmetries
\begin{equation}\label{super_string_sym}
c_a^{\vL,\sigma}=c_{-a}^{\vL,\sigma}=c_{a+2N}^{\vL,\sigma}
=c_{N-a}^{\vL,1-\sigma}=c_a^{\vL',\sigma'},
\end{equation}
with $L'_j=L_{N-j}$ and $\sigma'\equiv\sigma+L_1+\cdots+L_{N-1} \pmod{2}$.
For $\vL=\ve_{\ell}$ ($0\leq \ell<N)$ it reduces to the 
level-$N$ A$^{(1)}_1$ string functions $c_a^{\ell}$~\cite{JM84,KP84},
\begin{equation}\label{string}
c_a^{\ell}=c^{\ve_{\l},\sigma}_{a+N\sigma}.
\end{equation}

\begin{lemma}\label{lem_limit2}
Let $N$ be a positive integer, $\vL\in\Integer_+^N$ and $|q|<1$. Define
$a,h,k_i~(1\leq i\leq h), \sigma_{k_i}$ and $K,\bar{K}$ as in 
definition~\ref{def_b}. Then
\begin{equation}\label{limit2}
b_a^{\{L_k|k\in\bar{K}\}\{\sigma_k|k\in K\}}(q)=
\lim_{\substack{L_{k_1},\dots,L_{k_h}\to\infty\\ L_{k_i}\equiv
\sigma_{k_i}\pmod{2},~(1\leq i\leq h)}} T\bigl(\vL,\frac{a}{2}\bigr).
\end{equation}
\end{lemma}

\begin{proof}
This follows immediately from lemma~\ref{lem_explicit} noting that $x_j$ 
in \eqref{xjdef} may be rewritten (by inserting the explicit form of $m_0$ 
and replacing $a\to a/2$) as \eqref{xj} for all $1\leq j\leq N$. 
This implies for $|q|<1$,
\begin{equation}
\lim_{L_k\to\infty} x_j=\begin{cases} 0 & \text{for $j<k$}\\
x_j & \text{for $j\geq k$}. \tag*{$\qed$} \end{cases} 
\end{equation}
\renewcommand{\qed}{}
\end{proof}

Recalling that $x_N=1$, equation~\eqref{b_def} with $k_h=N$ simplifies to
\begin{equation}
b_a^{\{L_k|k\in\bar{K}\}\{\sigma_k|k\in K\}}(q)
=\frac{1}{(q)_{\infty}^h}
\sum_{\substack{m_i\in\Integer_+-\frac{1}{2}L_i,~i\in \bar{K}\\
m_i\in \Integer-\frac{1}{2}\sigma_i,~i\in K-\{N\} \\
\frac{a}{2N}+(C^{-1}\vm)_{N-1} \in \Integer+\frac{1}{2}\sigma_N}}
q^{\vm C^{-1} \vm} \prod_{j\in \bar{K}}\frac{1}{(q)_{\frac{L_j}{2}+m_j}}.
\end{equation}
When further $h=1$ we may shift $m_j\to m_j-L_j/2~(1\leq j<N)$
and obtain for $\vL\in\Integer_+^{N-1}$
\begin{equation}\label{b_c_rel}
b_a^{\{L_k\}\{\sigma\}}(q)=
q^{\frac{N}{4(N+2)}\vL C^{-1}\vL} c_a^{\vL,\sigma}(q).
\end{equation}
This leads to
\begin{corollary}\label{cor_super_string}
Let $N$ be a positive integer, $\vL\in\Integer_+^{N-1}$, $a\in\Integer$
and $\sigma=0,1$, such that $a+N(C^{-1}\vL)_{N-1}+N\sigma$ is even. Then
\begin{equation}\label{super_string_lim}
c_a^{\vL,\sigma}(q)=q^{-\frac{N}{4(N+2)}\vL C^{-1}\vL}
\lim_{\substack{L_N\to\infty\\ L_N\equiv\sigma \pmod{2}}}
T\left( (\vL,L_N),\frac{a}{2}\right).
\end{equation}
\end{corollary}

Finally, for $\vL=\ve_{\ell}~(1\leq \ell\leq N)$, the limit of the 
$q$-supernomials is related to the limit of the $q$-multinomials 
$T^{(N)}_{\ell}(L,a)$ of ref.~\cite{S96}, 
\begin{equation} \label{multi_limit}
\lim_{\substack{L\to\infty\\ L\equiv \sigma \pmod{2}}} T(L\ve_N+\ve_{\ell},a)
=q^{\frac{\ell(N-\ell)}{4N}}
\lim_{\substack{L\to\infty\\ L\equiv \sigma \pmod{2}}} 
T^{(N)}_{\ell}\bigl(L,\frac{\ell}{2}-a\bigr).
\end{equation}

\section{$q$-Supernomial coefficients and partitions}\label{sec3}
In the following a partition theoretic interpretation of the $q$-supernomials 
is provided. As we will see, this interpretation can be used to prove the 
recurrences of lemma~\ref{lem_rec} and the limits of lemma~\ref{lem_limit1} 
combinatorially. 
Before delving into the combinatorics of supernomials,
let us introduce a slight modification of the $q$-supernomial $T(\vL,a)$,
(compare with \eqref{duality}),
\begin{equation}
\tilde{T}(\vL,a)=q^{a\ell_1}\qbin{\vL}{a-\frac{\ell_N}{2}}_{1/q}.
\end{equation}
The explicit form for $\tilde{T}$ is easily found to be
\begin{equation}\label{qsupernomial3}
\tilde{T}(\vL,a)=\sum_{j_1+\cdots+j_N=a}
q^{\sum_{k=1}^N j_k (j_k+L_1+\cdots+L_{k-1})}
\qbin{L_N}{j_N}\qbin{L_{N-1}+j_N}{j_{N-1}}\cdots\qbin{L_1+j_2}{j_1}.
\end{equation}

\subsection{$q$-Supernomials and Durfee dissections}
\begin{figure}
\epsfxsize = 15cm
\centerline{\epsffile{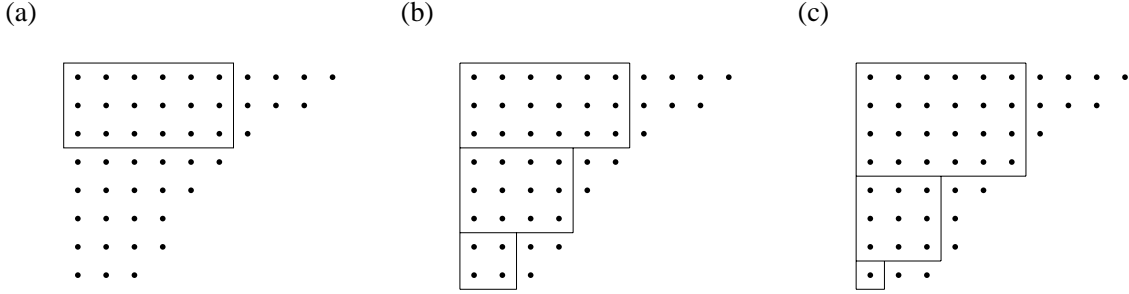}}
\caption{(a) Durfee rectangle of excess $3$ of the partition $\pi$.
(b) The $(3,1,0)$-Durfee dissection of $\pi$. (c)
The $(2,0,0)$-Durfee dissection of $\pi$.}\label{fig1}
\end{figure}

In order to interpret equation~\eqref{qsupernomial3} we
first give some basic definitions about Durfee dissections, 
generalizing some of the concepts introduced in refs.~\cite{Andrews79,W96}.
In the following a partition and its Ferrers graph are identified.
\begin{definition}
The Durfee rectangle of excess $E$ of a partition
is the maximal rectangle of nodes with
$E$ more columns than rows.
\end{definition}
As an example, the Durfee rectangle of excess $3$
of the partition $\pi=10+9+7+6+5+4+4+3$ is drawn in figure~\ref{fig1}(a).

For brevity it is often convenient to suppress the phrase
``of excess $E$'' and in the following we sometimes refer to
just ``the Durfee rectangle'', assuming that its excess
has been fixed.
\begin{definition}
The width of a Durfee rectangle is the number of columns.
The height of a Durfee rectangle is the number of rows.
\end{definition}

The part of a partition below its Durfee rectangle
is again a partition (with parts less than or equal to the width of
the Durfee rectangle) and we can define
a second Durfee rectangle, then a third and so on.
\begin{definition}
The $(E_1,\dots,E_n)$-Durfee dissection of a partition is obtained
by successively drawing Durfee rectangles of excess 
$E_1,E_2,\dots,E_n$.
\end{definition}
The $(3,1,0)$-Durfee dissection of $\pi$ is shown in figure~\ref{fig1}(b).

As a final definition we have
\begin{definition}\label{defadm}
Let $N$ be a positive integer and set $\vL\in\Integer_+^N$.
A partition is $(\vL;a)$-admissible if it has 
\begin{enumerate}
\item
exactly $a$ parts,
\item
no parts exceeding $L_1+L_2+\cdots + L_N$,
\item
no parts below its 
$(L_{N-1}+\cdots+L_1,\dots,L_2+L_1,L_1,0)$-Durfee dissection.
\end{enumerate}
The set of $(\vL;a)$-admissible partitions is denoted as
$S_{\vL;a}$.
\end{definition}
Using 1, the third condition is of course equivalent to the condition 
that the heights of the successive Durfee rectangles add up to $a$.

For $L \geq 7$, the partition $\pi$ is $(1,2,L;8)$-admissible as well as 
$(0,2,L+1;8)$-admissible, see figure~\ref{fig1}(b) and (c), respectively.
Note that the set $S_{L\ve_N;a}$ corresponds to the set of partitions 
with at most $N$ successive Durfee squares, such that the number of parts 
is $a$ and no part exceeds $L$.

We are now ready for the first claim of this section.
\begin{lemma}
For $\vL\in\Integer^N_+$,
$\tilde{T}(\vL,a)$ is the generating function of $(\vL;a)$-admissible
partitions.
\end{lemma}
\begin{proof}
Let $\lambda$ be an $(\vL;a)$-admissible partition, and let $j_{\ell}$ 
be the height of the $\ell$-th successive Durfee rectangle $D_{\ell}$ 
counted from below $(\ell=1,\dots,N)$. The number of nodes of $D_{\ell}$
is $j_{\ell}(j_{\ell}+L_1+\cdots+L_{\ell-1})$. The part of $\lambda$ to 
the right of $D_{\ell}$ (and below $D_{\ell+1}$ when $\ell<N$) is a 
partition with largest part $\leq j_{\ell+1}-j_{\ell}+L_{\ell}$ 
(with $j_{N+1}=0$)
and number of parts $\leq j_{\ell}$. Since the generating function of 
such partitions is given by $\qbin{L_{\ell}+j_{\ell+1}}{j_{\ell}}$ we 
find that the generating function of $(\vL;a)$-admissible partitions reads
\begin{equation*}
\sum_{j_1+\cdots + j_N=a} \prod_{\ell=1}^N
q^{j_{\ell}(j_{\ell}+L_1+\cdots + L_{\ell-1})}
\qbin{L_{\ell}+j_{\ell+1}}{j_{\ell}} = \tilde{T}(\vL,a).
\end{equation*}
Here the restriction in the sum over $j_{\ell}$ arises from the
fact that the sum over the respective heights of the Durfee  
rectangles should add up to $a$, to ensure that condition~3 of
definition~\ref{defadm} is satisfied.
\end{proof}

Having established the above result it is convenient to introduce a 
graphical notation for the set $S_{\vL;a}$. We proceed to do so in 
several steps. First, a rectangle
\begin{equation*}
\epsfysize = 1.5cm
\centerline{\epsffile{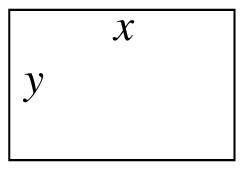}}
\end{equation*}
represents the partition of $xy$ consisting of $y$ rows of $x$ nodes. 
Alternatively, we can view the above rectangle as the (trivial) set of 
partitions with $y$ parts, all parts having size $x$.  Second, the triangle
\begin{equation*}
\epsfysize = 1.5cm
\centerline{\epsffile{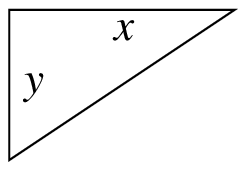}}
\end{equation*}
represents the set of partitions with at most $y$ parts, no part exceeding $x$.
We now combine these two elements to represent more complicated sets of 
partitions such as $S_{\vL;a}$.
As an example we first take $\vL=L\ve_1$, in which case $S_{\vL;a}$
corresponds to the set of partitions with exactly $a$ parts,
no part being less than $a$ and no part exceeding $L$.
Graphically this becomes
\begin{equation*}
S_{L\ve_1;a} = \raisebox{-0.65cm}{\epsfysize = 1.5cm \epsffile{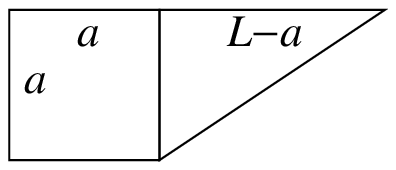}}
\end{equation*}
We assume that through these examples the idea is clear,
and for general $\vL$ we obtain the following
graphical representation of $S_{\vL;a}$: 
\begin{equation}\label{figsu}
S_{\vL;a}= \; \raisebox{-2.5cm}{\epsfysize=5cm \epsffile{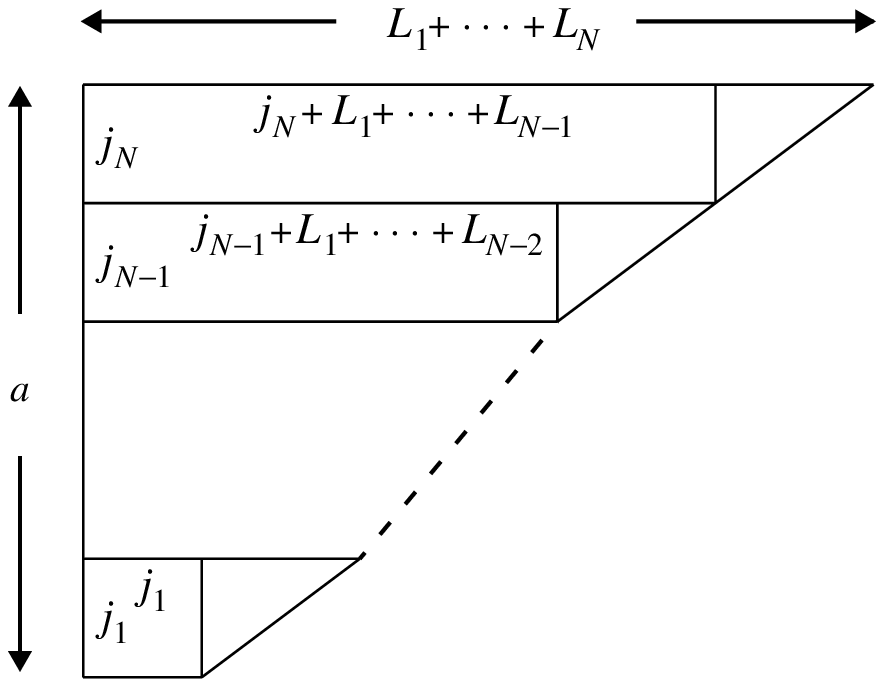}}
\end{equation}
To reiterate, to find a partition (or its Ferrers graph) 
in the set $S_{\vL;a}$ one interprets the above graph as follows.
Choose integers $j_1,\dots,j_N$ such
that their sum is $a$ and draw $N$ successive rectangles of 
$j_k$ times $j_k+L_1+\cdots+L_{k-1}$ nodes. Then, to the right of
the $k$-th rectangle from below, one may draw the Ferrers graph of
any partition that has no more than $j_k$ parts and that has largest
part not exceeding $L_k+j_{k+1}-j_k$ (with $j_{N+1}=0$). 

In the next section we make frequent use of (parts of) this graph.
When we refer to ``the $n$-th rectangle (triangle)'' this
is meant to indicate the $n$-th rectangle (triangle) counted from below. 
Also, to allow for a more compact notation, we often omit (part of)
the labels assuming that the excess of the $n$-th rectangle
is $L_1+\cdots + L_{n-1}$ unless stated otherwise.

\subsection{Combinatorial proof of lemma~\ref{lem_rec}}
The previously introduced graphical notation will aid us in asserting
lemma~\ref{lem_rec} purely combinatorially for $\vL-2\ve_n\in\Integer_+^N$.
In terms of $\tilde{T}(\vL,a)$ the recurrences \eqref{eqn_Trec} read
\begin{align}\label{qrec1} 
\tilde{T}(\vL,a) & =\tilde{T}(\vL+\ve_{n-1}-2\ve_n+\ve_{n+1},a)
+ q^{2a-n+\sum_{i=1}^{n-1}(n-i)L_i}\tilde{T}(\vL-2\ve_n,a-n), \\
\intertext{for $n=2,\dots,N-1$, and}
\label{qrec2}
\tilde{T}(\vL,a) & =q^a \, \tilde{T}(\vL-2\ve_1+\ve_2,a)
+q^{2a-1}\tilde{T}(\vL-2\ve_1,a-1).
\end{align}

\subsubsection{Recurrences \eqref{qrec1}}
Inspection of equation \eqref{qrec1} shows that on the left-hand side we
have the generating function of $(\vL;a)$-admissible partitions, whereas
on the right-hand side we have the generating function of
$(\vL+\ve_{n-1}-2\ve_n+\ve_{n+1};a)$-admissible partitions
plus the generating function of $(\vL-2\ve_n;a-n)$-admissible
partitions multiplied by $q$ to the power $2a-n+\sum_{i=1}^{n-1}(n-i)L_i$.
We now show that 
$S_{\vL+\ve_{n-1}-2\ve_n+\ve_{n+1};a}\subseteq S_{\vL;a}$ and,
setting $T_{\vL;a}=
S_{\vL;a}-S_{\vL+\ve_{n-1}-2\ve_n+\ve_{n+1};a}$,
that there is a bijection between
$T_{\vL;a}$ and $S_{\vL-2\ve_n;a-n}$, such that 
we have to remove $2a-n+\sum_{i=1}^{n-1}(n-i)L_i$ nodes
from each element of $T_{\vL;a}$ 
to map it to an element of $S_{\vL-2\ve_n;a-n}$.

To show that $S_{\vL+\ve_{n-1}-2\ve_n+\ve_{n+1};a}\subseteq S_{\vL;a}$
we make repeated use of the identities
\begin{align} &
\raisebox{-0.8cm}{\epsfysize=1.6cm \epsffile{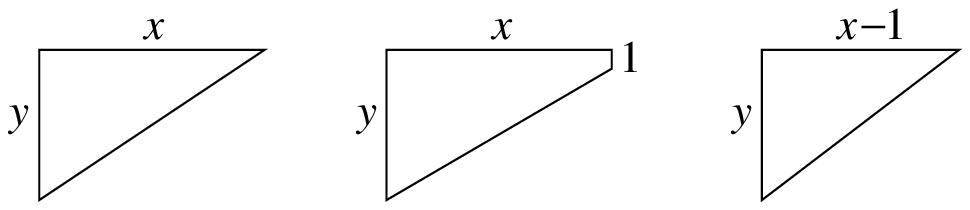}}
\hspace{-2.4cm} \dot{\cup}
\hspace{-3.5cm} = &
\label{t1} \\
\intertext{and}
& \raisebox{-0.8cm}{\epsfysize=1.8cm \epsffile{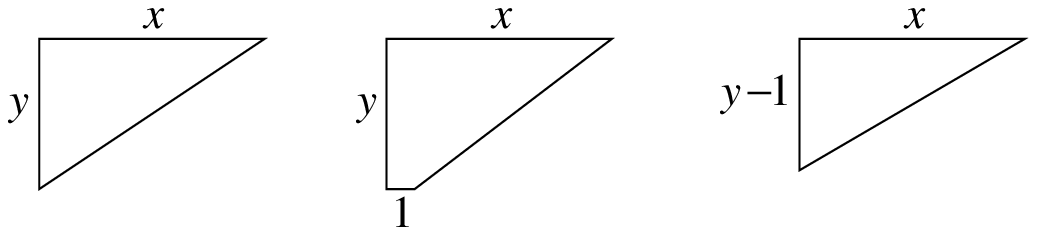}}
\hspace{-3.1cm} \dot{\cup}
\hspace{-3.2cm} = & \label{t2}
\end{align}
where it will be customary not to explicitly draw the rectangles of size
$1\times y$ and $x\times 1$. The symbol $\dot{\cup}$ stands for the disjoint
union.  Of course, \eqref{t1} just reflects the fact
that a partition with at most $y$ parts and  no part exceeding $x$
is either a partition with largest part precisely $x$,
or a partition with no parts exceeding $x-1$ (with in both cases no more
than $y$ parts),
and translated into generating functions \eqref{t1}  corresponds to
\eqref{bi1} with $L=x+y$ and $a=x$.
Similarly, \eqref{t2} encodes that a partition with at most
$y$ parts and no part exceeding $x$ is either a partition 
with precisely $y$ parts or a partition with at most $y-1$ parts
(with in both cases no part exceeding $x$),
corresponding to \eqref{bi2} with $L=x+y$ and $a=x$.

Using \eqref{t1} we may write
\begin{equation*}
\eqref{figsu} =
\raisebox{-0.7cm}{\epsfxsize=4.3cm\epsffile{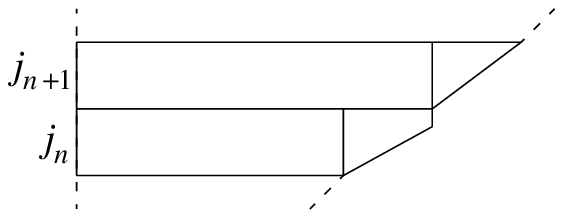}} \dot{\cup}\; \:
\raisebox{-0.7cm}{\epsfxsize=4.3cm\epsffile{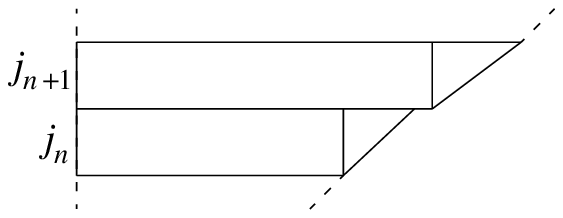}}
\end{equation*}
We now deform the first graph by increasing the height of the 
$(n+1)$-th rectangle by one thereby decreasing the height of
$n$-th rectangle. Since this only changes the interior
of the graph, it still represents the same set of partitions.
After this change, we shift variables $j_{n+1}\to j_{n+1}-1$
and $j_n\to j_n+1$, so that $j_n$ and $j_{n+1}$ again label
the heights of the $n$-th and $(n+1)$-th rectangle.
The second graph is rewritten using \eqref{t2} on the $n$-th
triangle.  In all this gives that \eqref{figsu} equals
\begin{equation*}
\raisebox{-0.7cm}{\epsfxsize=4.3cm \epsffile{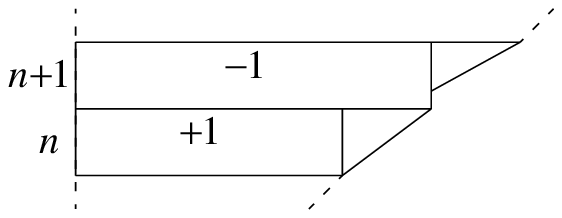}} \dot{\cup} \; \:
\raisebox{-0.7cm}{\epsfxsize=4.3cm \epsffile{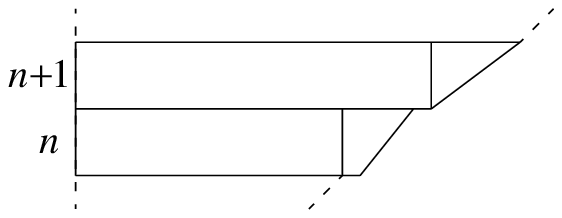}} \dot{\cup} \; \:
\raisebox{-0.7cm}{\epsfxsize=4.3cm \epsffile{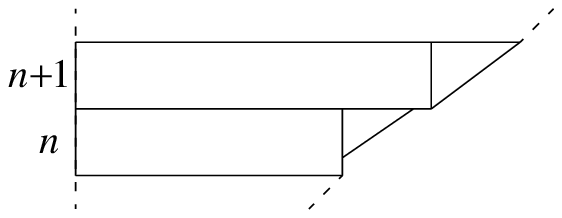}}
\end{equation*}
omitting the labels $j_n$ and $j_{n+1}$. Instead we use just the label
$k$ for the $k$-th rectangle (and triangle). The $-1$ and $+1$ in the first
graph are to indicate that the excesses of the $(n+1)$-th and $n$-th
rectangles are $L_1+\cdots +L_n-1$ and $L_1+\cdots +L_{n-1}+1$, respectively,
deviating $-1$ and $+1$ from their default values of equation \eqref{figsu}.

Now the second graph is deformed by increasing the width of the
$n$-th rectangle and decreasing the width of the
$(n+1)$-th rectangle, both by one unit.  
The third graph is rewritten using \eqref{t2} on the $(n-1)$-th
triangle.
As a result 
\begin{equation*}
\begin{split}
\eqref{figsu}= \; &
\raisebox{-0.7cm}{\epsfxsize=4.3cm \epsffile{fig10.ps}}  \dot{\cup}\;
\raisebox{-1cm}{\epsfxsize=4.3cm \epsffile{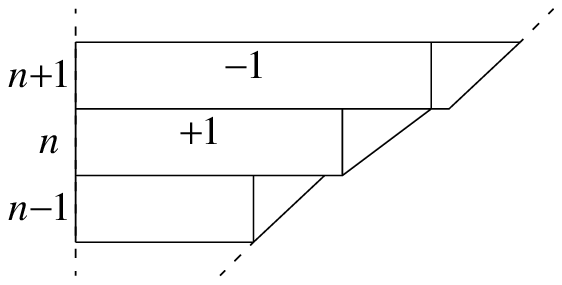}}  \\
\dot{\cup}  \; & 
\raisebox{-1cm}{\epsfxsize=4.3cm \epsffile{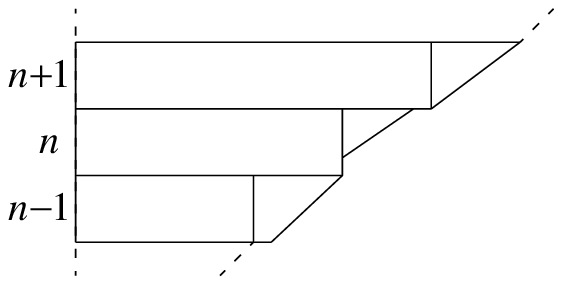}} \dot{\cup}\;
\raisebox{-1cm}{\epsfxsize=4.3cm \epsffile{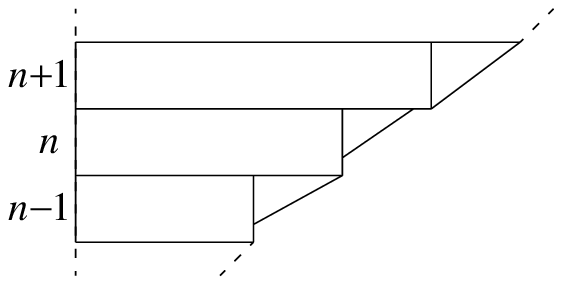}}
\end{split}
\end{equation*}
Now the third graph is deformed by increasing the height and width of
the $(n-1)$-th rectangle by one, such that the height of the $n$-the
rectangle is decreased by one. Further we decrease the width of the
$(n+1)$-th rectangle.
Finally, the last graph is rewritten using \eqref{t2} on the $(n-2)$-th
triangle.

This is the general pattern and we gradually move down in the graphs.
That is, after using \eqref{t2} on the $k$-th triangle of 
the last ($=(n-k+2)$-th)
graph, we get two new graphs. The first is deformed by increasing the
height and width of the $k$-th rectangle by one by shifting the rectangles
$k+1$ to $n-1$ one unit up and by decreasing the height of the $n$-th rectangle
by one. Furthermore the width of the $(n+1)$-th rectangle is decreased by one.
The second graph is again split using \eqref{t2} on the $(k-1)$-th triangle.
Iterating  this process we end up with $n+2$ graphs.
The first $n+1$ (which are the graphs that have all been deformed) read
\begin{equation*}
\begin{split}
&
\raisebox{-0.6cm}{\epsfxsize=4.3cm \epsffile{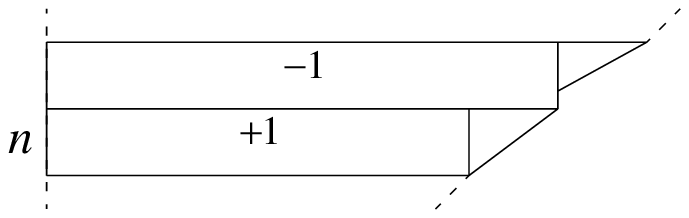}} \dot{\cup}
\raisebox{-0.8cm}{\epsfxsize=4.3cm \epsffile{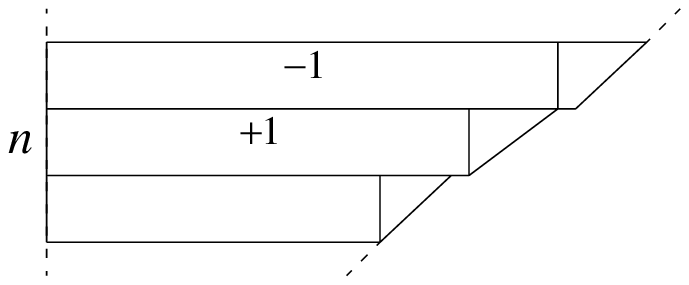}} \dot{\cup}
\raisebox{-1cm}{\epsfxsize=4.3cm \epsffile{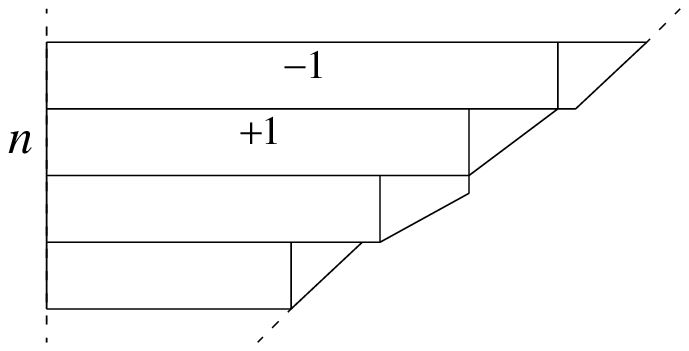}} \\
& \qquad  \dot{\cup} \cdots \dot{\cup}
\raisebox{-1.2cm}{\epsfxsize=4.3cm \epsffile{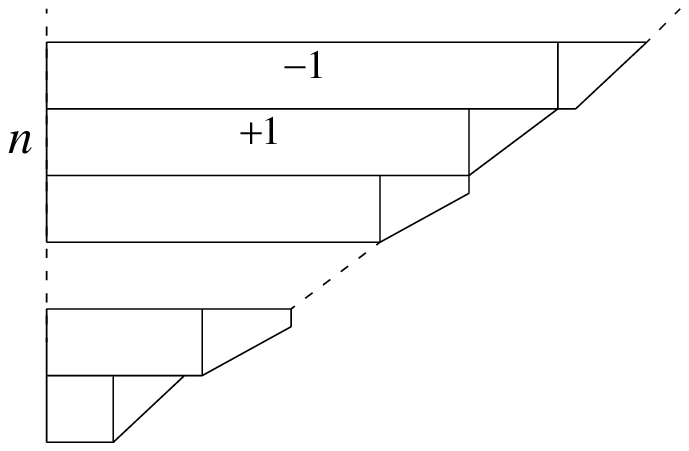}} \dot{\cup}
\raisebox{-1.2cm}{\epsfxsize=4.3cm \epsffile{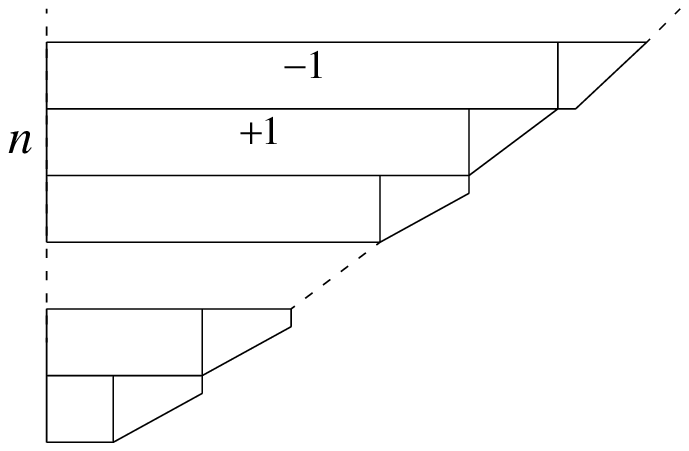}}
\end{split}
\end{equation*}
Now note that the last two graphs can be combined
by application of \eqref{t1} on the first triangles.
After this we can again combine the last two graphs,
now by use of \eqref{t1} on the second triangles.
Continuing this all the way up, all graphs
are combined into a single graph, by $(n-1)$-fold application
of \eqref{t1} followed by the single use of \eqref{t2}
on the $(n+1)$-th triangle. The result of all this is
the graph
\begin{equation}\label{setS}
\raisebox{-0.5cm}{\epsfxsize=4cm \epsffile{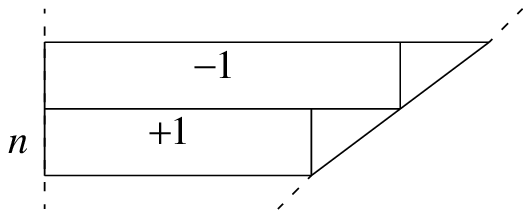}}
\end{equation}
which represents precisely the set of 
$(\vL+\ve_{n-1}-2\ve_n+\ve_{n+1};a)$-admissible partitions.

Let us now consider the neglected $(n+2)$-th graph representing
the set $T_{\vL;a}=S_{\vL;a}-S_{\vL+\ve_{n-1}-2\ve_n+\ve_{n+1};a}$.
The non-trivial part of this graph reads
\begin{equation}\label{setT}
\raisebox{-1.2cm}{\epsfxsize=4.3cm \epsffile{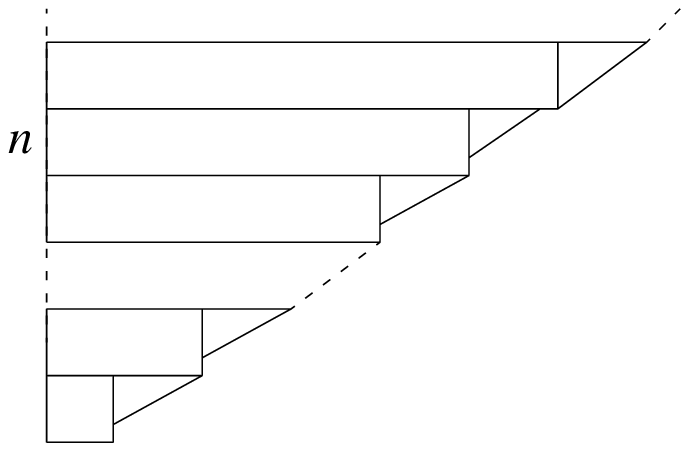}}
\end{equation}
We wish to show that there is a bijection between the set represented by
this graph and the set $S_{\vL-2\ve_n;a-n}$, such that each element of 
$T_{\vL;a}$ is mapped to an element of $S_{\vL-2\ve_n;a-n}$ by removal of 
$2a-n+\sum_{i=1}^{n-1}(n-i)L_i$ nodes.
In fact, it turns out that the map from a partition in $T_{\vL;a}$ to a 
partition in $S_{\vL-2\ve_n;a-n}$ acts only on the Durfee rectangles, and is
independent of the particular details of the partition.
We therefore temporarily change our interpretation of the graphs, and 
view the graph in equation \eqref{setT} as not representing
the set $T_{\vL;a}$, but as representing an arbitrary element of $T_{\vL;a}$.
That is, at the position of the $k$-th the triangle we assume an arbitrary
partition with number of parts less than or equal to the height of the 
triangle and
largest part less than or equal to the width of the triangle.
Then we remove the last row of the the first $n$ rectangles,
the first column and the column just preceding the $(n+1)$-th triangle.
Marking these rows and columns in grey, we get
\begin{equation*}
\raisebox{-1.5cm}{\epsfxsize=5cm \epsffile{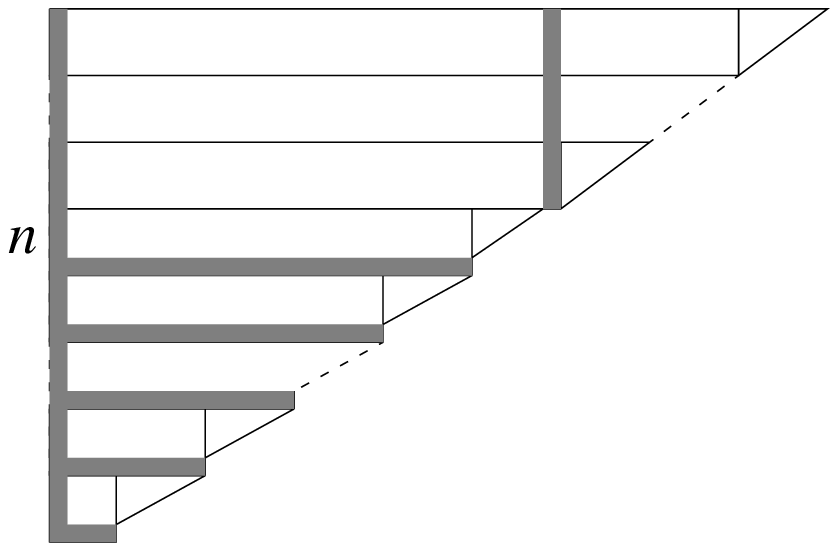}}
 \to \raisebox{-1.5cm}{\epsfxsize=5cm \epsffile{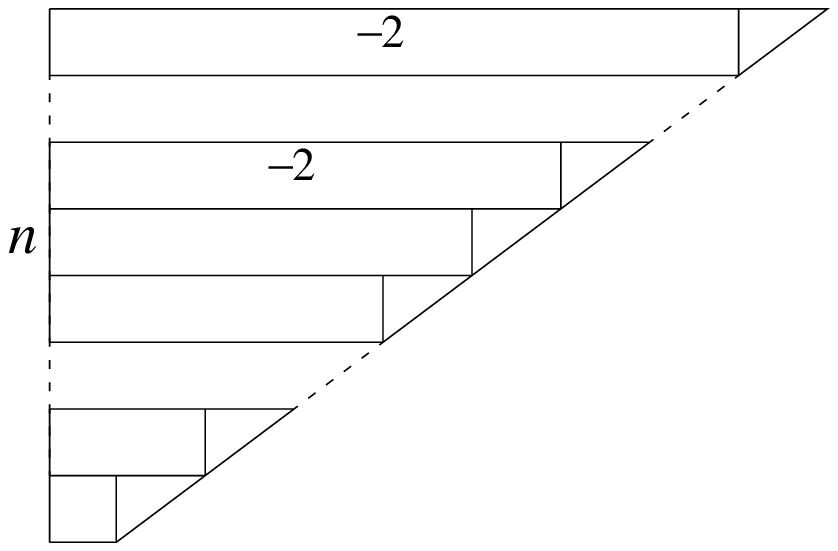}}
\end{equation*}
Since we have removed $n$ rows and two columns the right-hand side has
$a-n$ rows and no more than $L_1+\cdots +L_N-2$ columns,
and we conclude that it represents the set $S_{\vL-2\ve_n;a-n}$ (or arbitrary 
partitions in that set, depending which interpretation of the graph is chosen).
Since our map is trivially reversible, we are done.

\subsubsection{Recurrence \eqref{qrec2}}
The proof of the exceptional case \eqref{qrec2} proceeds in identical fashion
to the proof of \eqref{qrec1}.
The only modification arises from the fact that the set $S$, represented
by the graph of equation \eqref{setS} with $n=1$ 
(so that the dotted lines below the graph should be omitted),
does not represent the set $S_{\vL-2\ve_1+\ve_2;a}$. 
However, the map onto
$S_{\vL-2\ve_1+\ve_2;a}$ is trival, since we only have to remove the
first column of \eqref{setS}. 
Since this corresponds to removal of $a$ nodes from each element of $S$
this gives rise to the additional factor $q^a$ in the first term
of the right-hand side of \eqref{qrec2}.
Obviously this map is reversible, completing the proof.

\subsubsection{Recurrence for $n=N$}
The recurrences in lemma~\ref{lem_rec} hold for $n=1,\dots,N-1$.
Using the above graphical methods, it is an easy matter to verify that
on setting $n=N$ one gets the modified result
\begin{align}\label{ex1}
\tilde{T}(\vL,a)
&= q^{\ell_1}\tilde{T}(\vL+\ve_{N-1}-\ve_N,a-1)
+ q^{2a-N+\sum_{i=1}^N(N-i)L_i}\tilde{T}(\vL-2\ve_N,a-N) \notag \\
& \qquad + \tilde{T}(\vL+\ve_{N-1}-2\ve_N,a).
\end{align}
Similarly one can find that
\begin{align}\label{ex2}
\tilde{T}(\vL,a)
&= \tilde{T}(\vL+\ve_{N-1}-\ve_N,a)
+ q^{2a-N+\sum_{i=1}^N(N-i)L_i}\tilde{T}(\vL-2\ve_N,a-N) \notag \\
& \qquad 
+ q^{2a-N-1+\sum_{i=1}^N(N-i+1)L_i}\tilde{T}(\vL+\ve_{N-1}-2\ve_N,a-N-1).
\end{align}
Both these results, which hold for 
$N\geq 2$ and $\vL-2\ve_N \in \Integer^N_+$, will be briefly touched upon 
in section~\ref{sec:42}. However, contrary to the cases $1 \leq n \leq N-1$, 
neither will be essential in any of the proofs of subsequent theorems.

\subsection{A Durfee rectangle like identity}
A last combinatorial aspect of the $q$-supernomials we wish to mention
is that one can prove lemma~\ref{lem_limit1} by interpreting equation
\eqref{limit1} as a successive Durfee rectangle identity.

Let us first recall the Durfee rectangle identity
\begin{equation*}
\sum_{n=m}^{\infty} \frac{q^{(n-m)(n+m)}}{(q)_{n-m}(q)_{n+m}}
=\frac{1}{(q)_{\infty}},
\end{equation*}
true for all (integer) $m\geq 0$.
This identity can be understood by noting that the summand
is the generating function of partitions with Durfee rectangle
of width $n+m$ and height $n-m$.
Hence summing over $n\geq m$ yields the generating function of all
partitions.

Returning to the $q$-supernomial of equation~\eqref{qsupernomial1}
and to lemma~\ref{lem_limit1}, we note that 
\begin{equation*}
q^{\sum_{k=2}^N j_{k-1} (L_k+\cdots+L_N-j_k)}
\qbin{L_N}{j_N}\qbin{L_{N-1}+j_N}{j_{N-1}}\dots\qbin{L_1+j_2}{j_1}
\end{equation*}
is the generating function of partitions with at most $j_1+\cdots+j_N$
parts, no part exceeding $L_1+\cdots+L_N-j_1$, with $N-1$
successive Durfee rectangles of heights $j_1,\dots,j_{N-1}$ and widths 
$L_2+\cdots+L_N-j_2,L_3+\cdots+L_N-j_3,\dots,L_N-j_N$.
If we denote the set of such partitions by $S_{\vL;\bs j}$
we get, using our graphical notation for restricted partition sets,
\begin{equation*}
S_{\vL;\bs j}= \; \raisebox{-2.5cm}{\epsfysize=5cm \epsffile{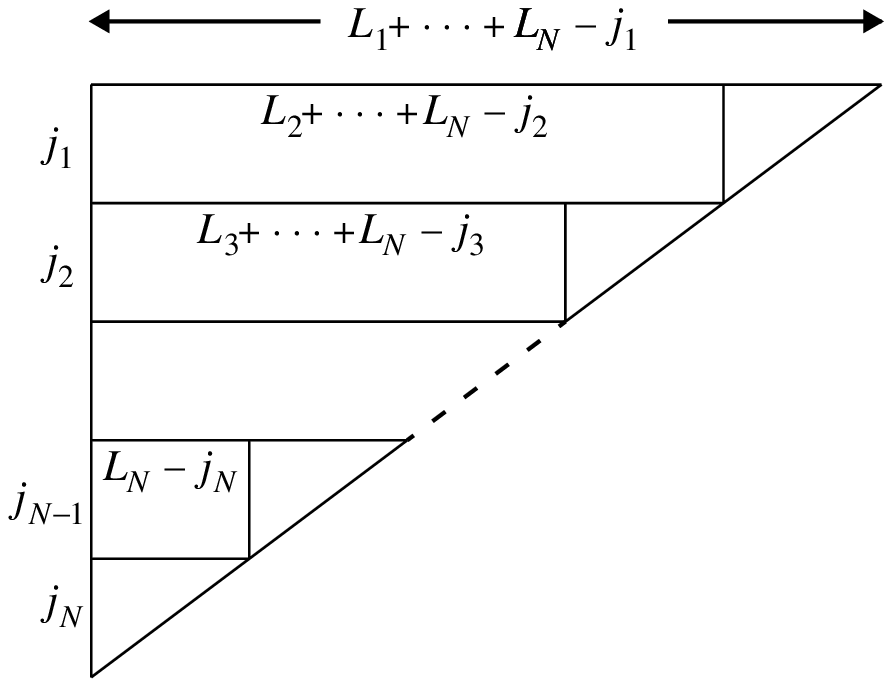}}
\end{equation*}
Sending $L_m\to\infty$, the generating function of $S_{\vL;\bs j}$ is
non-vanishing only when $j_1=\dots=j_{m-1}=0$.
But in that case $\lim_{L_m\to\infty}
S_{\vL;\bs j}$ becomes the set of partitions with at most $j_m+\cdots+j_N$
parts, $N-m$ successive Durfee rectangles of heights $j_m,\dots,j_{N-1}$ and 
widths $L_{m+1}+\cdots+L_N-j_{m+1},\dots,L_N-j_N$. Hence the set
$$
\lim_{L_m\to\infty} 
\underset{\substack{j_m+\cdots+j_N=a \\ j_1,\dots,j_{m-1}=0}}{\bigcup}
S_{\vL;\bs j}
$$
corresponds to the set of all partitions with at most $a$ parts, so that
\begin{equation*}
\lim_{L_m\to\infty} \qbin{\vL}{a-\frac{\ell_N}{2}}=\frac{1}{(q)_a}.
\end{equation*}

\section{The Andrews--Gordon identities and $q$-supernomials}\label{sec4}
We now come to the second part of our paper where it is shown that
several well-known $q$-series identities have a natural extension
using the $q$-supernomials.
The first and simplest example is Andrews' analytic
counterpart of Gordon's partition theorem~\cite{Andrews74}, for which we
give polynomial analogues depending on multiple finitization parameters.

\subsection{Brief history}\label{sec41}
To put things in a somewhat broader perspective, we first recall
some facts about the Rogers--Ramanujan identities.
For $a=0,1$ define the function
\begin{equation}\label{fdef}
f(x;q) = f(x) =  \sum_{n\in\Integer} \frac{q^{n(n+a)} x^{2n}}{(x)_{n+1}} .
\end{equation}
Obviously, 
\begin{equation*}
\lim_{x\to 1^{-}} (1-x) f(x;q) = 
\sum_{n\geq 0} \frac{q^{n(n+a)}}{(q)_n},
\end{equation*}
which equals one side of the Rogers--Ramanujan identities.
Now view $f$ as the generating function for polynomials, setting
\begin{equation*}
f(x) = \sum_{M\in\Integer} P(M) x^{M} .
\end{equation*}
By the $q$-binomial theorem~\cite{Andrews76}
\begin{equation}\label{qbinthm}
\frac{1}{(x)_{n+1}}=\sum_{j\in\Integer}\qbin{j+n}{j}x^j
\end{equation}
it is readily seen that
\begin{equation}\label{MM}
P(M) = \sum_{n\in\Integer} q^{n(n+a)} \qbin{M-n}{M-2n}.
\end{equation}
For $M\geq 0$ these polynomials were introduced by MacMahon~\cite{MacMahon} 
and in this case the right-hand side of \eqref{MM} corresponds to
the generating function of 
partitions with difference between parts at least two,
no parts exceeding $M+a-1$ and all parts exceeding $a$.

An altogether different representation for $P$ was found 
by Schur~\cite{Schur17},
\begin{equation}\label{S}
P(M)=\sum_{j=-\infty}^{\infty}
(-1)^j q^{(5j+2a+1)j/2} \qbin{M+a}{\lfloor \frac{1}{2}(M-5j) \rfloor},
\end{equation}
with $\lfloor . \rfloor$ the integer part function. Equating~\eqref{MM} 
and~\eqref{S} yields polynomial identities due to Andrews~\cite{Andrews70}. 
Letting $M\to \infty$ (for $|q|<1$), the Rogers--Ramanujan identities
are reproduced.

Apart from the fact that the function $f$ in \eqref{fdef}
acts as the generating function for the polynomials $P$,
it has some further interesting properties~\cite{Andrews85}.
Here we mention the $q$-difference equation
\begin{equation}\label{qdiff}
(1-x) f(x) = q^{a+1} x^2 f(qx)
\end{equation}
which implies the recurrences
\begin{equation*}
P(M) = P(M-1) + q^{M+a-1} P(M-2).
\end{equation*}
It is indeed these recurrences and the initial conditions $P(0)=P(1)=1$
which lead to the proof that \eqref{MM} and \eqref{S} are equal.

\subsection{The multivariable case}\label{sec:42}
In ref.~\cite{Andrews85} it was suggested that the following 
multivariable function might be a relevant generalization of $f$:
\begin{equation}\label{deff}
f(\vx) = \sum_{n_1,\dots,n_k \in\Integer} \frac{
q^{N_1^2+\cdots +N_k^2+N_a+\cdots+N_k}
x_1^{2N_1} x_2^{2N_1+2N_2} \cdots x_k^{2N_1+\cdots+2N_k}}
{(x_1)_{n_1+1} (x_2)_{n_2+1} \cdots (x_k)_{n_k+1}},
\end{equation}
with
\begin{equation*}
N_j = n_j+\cdots +n_k,
\end{equation*}
$\vx=(x_1,\dots,x_k)$ and $1 \leq a \leq k+1$.
For $k=2$ this function has been analysed in relation with a generalization
of the hard-hexagon model~\cite{AB87}. 
It is easy to show that $f$ obeys the $q$-difference equations
\begin{equation}\label{qdiff1}
(1-x_p) f(\vx) = q^{2p-\min\{a-1,p\}}
\Bigl( \prod_{j=1}^k x_j^{2\min\{j,p\}} \Bigr)
f(x_1,\dots,qx_p,\dots,x_k),
\end{equation}
for $p=1,\dots,k$. These equations generalize \eqref{qdiff}.

Whereas \eqref{fdef} is closely related to the 
Rogers--Ramanujan identities, \eqref{deff} has a connection to the 
Andrews--Gordon identities~\cite{Andrews74} given by
\begin{theorem}
For $k \geq 1$ and $1 \leq a \leq k+1$,
\begin{equation}\label{AG}
\sum_{n_1,\dots,n_k \geq 0} \! \! 
\frac{q^{N_1^2+\cdots+N_k^2+N_a+\cdots+N_k}}{(q)_{n_1} \cdots (q)_{n_k}}
=\frac{1}{(q)_{\infty}}\sum_{j=-\infty}^{\infty} (-1)^j q^{((2k+3)(j+1)-2a)j/2}
= \hspace{-4mm} 
\prod_{\substack{j=1 \\ j\not\equiv 0,\pm a \pmod{2k+1}}}^{\infty}
\hspace{-10mm} (1-q^j)^{-1}.
\end{equation}
\end{theorem}
Indeed we see that $\lim_{\vx\to(1^-,\dots,1^-)} (1-x_1)\dots (1-x_k) f(\vx)$
returns the left-hand side of \eqref{AG}.

Again viewing $f$ as a generating function for polynomials, we define $P$ as 
\begin{equation}\label{defP}
f(\vx) = \sum_{\vM \in \Integer^k} P(\vM) x_1^{M_1} \cdots x_k^{M_k}.
\end{equation}
Multiple application of the $q$-binomial theorem \eqref{qbinthm}
shows that for all $\vM \in \Integer^k$,
\begin{equation}\label{Pf}
P(\vM)=\sum_{n_1,\dots,n_k \in\Integer}
q^{N_1^2+\cdots+N_k^2+N_a+\cdots+N_k}
\prod_{j=1}^k \qbin{M_j+n_j-2\sum_{m=1}^j N_m}{M_j-2\sum_{m=1}^j N_m}.
\end{equation}
Less obvious is the following alternative form for $P$, generalizing
\eqref{S} and starring the $q$-supernomials.
\begin{theorem}\label{thmB}
Set $\vL=T\vM-\ve_{a-1}+\ve_k\in \Integer_+^k$ and $\vl=T^{-1}\vL$
with $T$ defined in section \ref{sec_prelims} with $N$ replaced by $k$.
Then for $k\geq 1$ and $1 \leq a \leq k+1$,
\begin{equation}\label{Pb}
P(\vM)=\sum_{j=-\infty}^{\infty} \biggl\{ q^{j((2k+3)(2j+1)-2a)}
\qbin{\vL}{\frac{b-a}{2}+(2k+3)j}
-q^{(2j+1)((2k+3)j+a)} \qbin{\vL}{\frac{b+a}{2}+(2k+3)j}\biggr\},
\end{equation}
with $b=k+1$ if $\ell_k+a+k$ is odd and $b=k+2$ if $\ell_k+a+k$ is even.
\end{theorem}

Before proving this theorem we note that equating the two
different representations for $P$ indeed yields a 
polynomial analogue of the Andrews--Gordon identity \eqref{AG}
based on $k$ finitization parameters.
Taking the limit $L_p\to\infty$ ($p=1,\ldots,k$) of $P(\vM)$
in the representation~\eqref{Pb} reduces to the middle form of~\eqref{AG}
by lemma~\ref{lem_limit1}. The same limit with $P(\vM)$ in representation
\eqref{Pf} reduces to the left-hand side of \eqref{AG}, recalling that
$\vM=T^{-1}(\vL+\ve_{a-1}-\ve_k)$.

We also remark that for $\vL=L \ve_1$, the polynomial identities implied 
by theorem~\ref{thmB} were proven by Foda and Quano~\cite{FQ95} and by 
Kirillov~\cite{Kirillov95}. (For $k=1$ these are
Andrews' identities discussed in section~\ref{sec41}.) 
In ref.~\cite{W96} (a slight variation of) the polynomial 
identities implied by the theorem were conjectured for all 
$\vL = L \ve_p$ ($p=2,\dots,k$), and proven for $\vL=L \ve_k$.

\begin{proof}[Proof of theorem \ref{thmB}]
We show that both \eqref{Pf} and \eqref{Pb} satisfy the recurrences
\begin{equation}\label{recP}
P(\vM) = P(\vM-\ve_p) + q^{{M_p}-\min\{a-1,p\}} P(\vM- 2 T^{-1} \ve_p),
\end{equation}
for $p=1,\dots,k-1$ and all $\vM\in\Integer^k$. The representation~\eqref{Pf}
satisfies this recurrence due to~\eqref{qdiff1} and~\eqref{defP}. To show
that~\eqref{Pb} obeys~\eqref{recP}, we set 
$P(\vM)=P(T^{-1}(\vL+\ve_{a-1}-\ve_k))=\tilde{P}(\vL)$ which implies
\begin{equation}\label{recPt}
\tilde{P}(\vL)  = \tilde{P}(\vL+\ve_{p-1}-2\ve_p+\ve_{p+1})
+ q^{{\ell_p}-p} \tilde{P}(\vL- 2\ve_p),
\end{equation}
for $\vL\in \Integer^k$, $\vl=T^{-1}\vL$ and $p=1,\dots,k-1$.
The proof that \eqref{Pb} satisfies \eqref{recPt} follows
immediately from \eqref{eqn_rec}.

Thanks to these recurrences and lemma~\ref{lem_complete} it suffices to show 
that equations~\eqref{Pf} and \eqref{Pb} are equal for the initial condition
$\vL=L\ve_1$ with $L\geq 0$.
This, however, is the one-variable polynomial identity
proven in refs.~\cite{FQ95,Kirillov95}.
\end{proof}

Having established theorem~\ref{thmB}, we wish to make a remark
about the proof. In order 
to establish the equivalence between the representations 
\eqref{Pf} and \eqref{Pb} for $P(\vM)$, we have used the recurrences
\eqref{recP} for $p=1,\ldots,k-1$ plus the one-parameter family of 
polynomial identities of refs.~\cite{FQ95,Kirillov95}. 
Of course this second ingredient
is a rather non-trivial result and one may wonder whether theorem~\ref{thmB}
can be obtained without it. Indeed we point out that both
\eqref{Pf} and \eqref{Pb} satisfy the recurrence \eqref{recP} 
when $p=k$. Since this means that we have $k$ independent recurrence
relations it now suffices to take only a finite number of initial
points to prove the theorem. However, admittedly, it is a non-trivial
task to determine a suitable set of initial points, due to the
``non-orthogonal'' nature of the recurrences. 
To prove that~\eqref{Pb} satisfies the $p=k$ recurrence, one needs
to employ the formulae \eqref{ex1} and \eqref{ex2} (or rather, their
corresponding $q\to 1/q$ analogue).
Treating the $b=k+1$ and $b=k+2$ cases separately, the proof
is a matter of straightforward algebra.

To conclude this section we mention that the function $P$ in
the representation of equation \eqref{Pf} has been extensively studied by
Bressoud in refs.~\cite{Bressoud89,Bressoud92} in his
work on the O'Hara--Zeilberger identity. One of the results
established in ref.~\cite{Bressoud89} is a partition theoretic
interpretation of $P(\vM)$ when $\vM\in\Integer_+^k$, obtained from 
theorem~2 therein by
dropping the restriction that the partitions should have $j$ parts
and by replacing $n_r$ by $M_r$ ($r=1,\dots,k$). 
(Note that in theorem~2 of ref.~\cite{Bressoud89} it should be added that
the number of parts is $j$ and $ns(B)$ should be corrected to $n_{s(B)}$.)
A second result from ref.~\cite{Bressoud89} is an interpretation
of $P(\vM)$ in terms of lattice paths, obtained from lemma~4 therein
by replacing $n_r$ by $M_r$ ($r=1,\ldots,k$) and
summing over $m_1,\dots,m_k \geq 0$.

\section{Supernomial boson--fermion identities and continued fractions}
\label{sec5}
\subsection{$(\vm,\vn)$-Systems}\label{sec_mn}
In the previous section it was demonstrated how the $q$-supernomials can be 
applied to obtain polynomial identities related to identities of the 
Rogers--Ramanujan-type.
In the remainder of this paper we treat a much more general application
of the $q$-supernomials, which takes the results of section~\ref{sec4} as
special case. In particular we show that to each pair of integers $p$ and
$k$, with $0<2k<p$ and $\gcd(p,k)=1$, 
there corresponds a family of polynomial identities depending on $N$ 
variables, with $N$ restricted by $0<N<(p-1)/k-1$.

Due to the complexity of this general result, we need a more systematic
way to express our formulae, and to motivate the notation
needed, we first recast some of our previous expressions in different
forms, which have their genesis in the recent literature on
``boson--fermion'' identities (see ref.~\cite{BMS96} and references
therein).

Returning to the representation \eqref{Pf} of the function $P$,
we replace $k$ by $N$ and
define an $N$-dimensional vector $\vn=(n_1,\dots,n_N)$ and a
``companion'' vector $\vm=(m_1,\ldots,m_N)$, such that
$m_j=M_j-2\sum_{\ell=1}^j N_{\ell}$.
Clearly, given $\vM$, the vector $\vm$ determines the vector $\vn$ and 
vice versa.
If we now eliminate the variables $N_j$ in the quadratic exponent of $q$,
equation \eqref{Pf} takes the form
\begin{equation*}
P(\vM) = \sum_{\vn\in \Integer^{N}} 
q^{\vn T^{-1} ( \vn + \ve_N - \ve_{a-1})} 
\prod_{j=1}^N \qbin{m_j+n_j}{m_j},
\end{equation*}
where $T$ is the Cartan-type matrix defined in section \ref{sec_prelims}.
Using $T \vM = \vL +\ve_{a-1}-\ve_N$, the relation between 
$\vn$ and $\vm$ can be rewritten in the canonical form of an 
``$(\vm,\vn)$-system''~\cite{Berkovich94},
\begin{equation*}
\vm+\vn = \frac{1}{2} ( \I_T \, \vm + \ve_{a-1} - \ve_N + \vL).
\end{equation*}
Here $\I_T$ is the incidence matrix of the tadpole graph with $N$ nodes,
$(\I_T)_{i,j}=\delta_{|i-j|,1}+\delta_{i,j}\delta_{i,N}$.

Now replacing $q\to 1/q$ and using equations~\eqref{dual_qbin}
and \eqref{duality}, the result of theorem~\ref{thmB}
translates to the following polynomial identity
\begin{multline}\label{dualAG}
q^{\frac{1}{4}(N-a+1)}
\sum_{\vn\in \Integer^N}
q^{\frac{1}{4}\vm T\vm-\frac{1}{2}(m_{a-1}-m_N)}
\prod_{j=1}^N \qbin{m_j+n_j}{m_j}\\
=q^{\frac{1}{4N}(b-a)^2}
\sum_{j=-\infty}^{\infty}\Bigl\{ q^{\frac{j}{N}(3pj+p(b-N)-3a)}
T\bigl(\vL,\tfrac{b-a}{2}+pj\bigr)
-q^{\frac{1}{N}(pj+a)(3j+b-N)}
T\bigl(\vL,\tfrac{b+a}{2}+pj\bigr) \Bigl\}, 
\end{multline}
for $\vL \in \Integer_+^N$,
$p=2N+3$, $a=1,\dots,N+1$ $(m_0=0$) and $b=N+1$ for $\ell_N+a+N$ odd
and $b=N+2$ for $\ell_N+a+N$ even, where $\vl=T^{-1}\vL$.

In the following we generalize the above polynomial
identity by considering more complicated incidence and Cartan-type matrices,
defined through the continued fraction expansion of $p/k$, the
simple case \eqref{dualAG} corresponding to $p=2N+3$ and $k=2$.
Before doing so, we first discuss identity \eqref{dualAG}, and its
subsequent generalizations, in the context of
solvable lattice models and conformal field theories.

\subsection{Conformal field theories, solvable lattice models and
boson--fermion identities}\label{sec52}
For each $\vL\in \Integer_+^N$, equation \eqref{dualAG} yields a 
polynomial identity. If we restrict $\vL$ to the lines  
$\vL = L \ve_m$ with $1\leq m \leq N$, the polynomials generated
by \eqref{dualAG} correspond to the one-dimensional
configuration sums of solvable lattices models of 
Date et al.~\cite{DJKMO87,DJKMO88} when their 
results are extended to ``non-physical'' regimes.
(For $N=1$ \eqref{dualAG} corresponds to the 4-state
Andrews--Baxter--Forrester model in regime IV~\cite{ABF84}.)

As was pointed out by Date et al., taking the limit $L\to\infty$
of the one-dimensional configuration sums yields branching functions
or characters of (coset) conformal field theories. 
The one-dimensional configuration sums are thus polynomial 
approximations of conformal characters, and are often referred to 
as ``finitized'' characters.
The $q$-supernomial identities involve
higher-dimensional $\vL$, providing a unification of 
various different ``finitized'' characters.
In fact, thanks to the multi-dimensionality of $\vL$, we can study much
more general $q$-series limits (by sending an arbitrary number of
components of $\vL$ to infinity), providing identities for
generalizations of the A$_1^{(1)}$ branching functions, 
see section~\ref{sec_limits}.

In the physics literature identities such
as \eqref{dualAG} are known as (polynomial) boson--fermion identities,
the left-hand side involving the $(\vm,\vn)$-system being ``fermionic'' and
the alternating-sign expression on the right-hand side being ``bosonic''.
Indeed, the left-hand side admits an interpretation as the partition function
for a system of quasiparticles with fractional statistics, obeying Fermi's
exclusion principle~\cite{KKMM93a,KKMM93b}.
The $(\vm,\vn)$-system of a fermionic function then
describes the combinatorics of the fermionic quasiparticles,
the $j$-th components of $\vn$ and $\vm$ having the interpretation of
occupation number of particles and antiparticles of charge $j$, 
respectively.

The term bosonic has its origin in conformal field theory (CFT).
Each CFT is described by an underlying symmetry algebra, such as the Virasoro 
algebra, the superconformal algebra or other extended algebras.
The states in a CFT can be determined from the highest
weight representation of the symmetry algebra by constructing the
Verma modules, each of which consists of a highest weight state and its
descendants, obtained by acting with the lowering operators on this state.
Some descendants may have zero norm and are referred to as
singular vectors. For the calculation of the characters or branching functions
the contribution of these singular vectors has to be subtracted out. 
This procedure is called the Feigin and Fuchs construction~\cite{FF82,FF83}
(and is analogous to the sieving technique of 
refs.~\cite{Andrews72,Andrews76}).
Since it relies on the use of a ``bosonic Fock space'',
the resulting character expressions are referred to as bosonic.
Even though the Feigin--Fuchs construction applies to conformal
characters and not their polynomial finitizations, the term bosonic is now 
generally used for expressions of alternating-sign-type,
such as the right side of \eqref{dualAG}.

When taking the finitization parameter(s) in polynomial
boson--fermion identities to infinity, one obtains $q$-series identities
of boson--fermion-type.
In certain special cases, the bosonic side can be further rewritten
using the triple product identity or its generalizations, to yield
identities of Rogers--Ramanujan type. The term boson--fermion-type is thus
referring to a wider class of identities than Rogers--Ramanujan-type.

We finally wish to emphasize that polynomial boson--fermion identities
which depend on a vector $\vL$ rather than a scalar $L$ have not
appeared in the literature before, and it is rather intriguing to
note that many of the existing polynomial identities related
to configuration sums of different solvable lattive models, can be
cast in one unifying form, using the $q$-supernomials.

\subsection{Continued fraction expansion and Takahashi--Suzuki decomposition}
\label{sec53}
After the previous intermezzo we describe how~\eqref{dualAG} can be 
generalized by considering the continued fraction expansion of $p/k$.
All notation and definitions of this subsection are borrowed from the work of 
Berkovich and McCoy~\cite{BM96}.

\begin{definition}\label{def_cf}
Let $p,k$ be integers such that $0<2k<p$ and $\gcd(p,k)=1$. 
Then the integers $n$ and $\nu_j$ $(0\leq j\leq n)$ are defined 
by the continued fraction expansion 
\begin{equation}\label{c_fraction}
\frac{p}{k}=1+ \nu_0+\cfrac{1}{\nu_1+\cfrac{1}{\nu_2+\ldots
+\cfrac{1}{\nu_n+2}}} \; .
\end{equation}
\end{definition}
When $k=1$, we have $n=0$ and the continued fraction 
expansion~\eqref{c_fraction} reads $p=\nu_0+3$.

Using $n$ and $\nu_j$ ($0\leq j\leq n$) we introduce the sums
\begin{equation}\label{t_m}
t_m=\sum_{j=0}^{m-1}\nu_j \quad (1\leq m\leq n+1)
\quad\text{and} \quad t_0=-1,
\end{equation}
which define a fractional incidence and Cartan-type matrix as follows.
\begin{definition}\label{def_IB}
Let $n$, $\nu_j$ $(0\leq j\leq n)$ and $t_m$ $(0\leq m \leq n+1)$ be 
given by the continued fraction expansion of $p/k$.
The fractional incidence matrix $\I_B$ is given by
\begin{equation} \label{I_B}
(\I_B)_{i,j} = \begin{cases}
\delta_{i,j+1} + \delta_{i,j-1} & \text{for
$1 \leq i<t_{n+1},~ i \neq t_m$}, \\
\delta_{i,j+1} + \delta_{i,j} - \delta_{i,j-1} & \text{for
$i=t_m,~1\leq m\leq n-\delta_{\nu_n,0}$}, \\
\delta_{i,j+1} + \delta_{\nu_n,0} \delta_{i,j} & \text{for $i=t_{n+1}$}.
\end{cases}
\end{equation}
$\I_B$ defines a Cartan-type matrix $B$ via $B=2 I -\I_B$,
with $I$ the $t_{n+1}$ by $t_{n+1}$ unit matrix.
\end{definition}
Notice that when $k=1$, the incidence matrix $\I_B$
has components $(\I_B)_{i,j}=\delta_{|i-j|,1}$ ($i,j=1,\dots,p-3$), so that
$B$ corresponds to the Cartan matrix of the Lie algebra A$_{p-3}$.

For $0\leq m\leq n$, we set the recurrences
\begin{align}\label{def_y}
y_{m+1}&=y_{m-1}+(\nu_m+2\delta_{m,n}+\delta_{m,0}) y_m, & y_{-1}&=0, &
y_0&=1,\\
\label{def_by}
\by_{m+1}&=\by_{m-1}+(\nu_m+2\delta_{m,n}+\delta_{m,0}) \by_m, & 
\by_{-1}&=-1, & \by_0&=1,
\end{align}
so that $y_{n+1}=p$ and $\by_{n+1}=p-k$.
This leads to the following definition.
\begin{definition}\label{def_taka}
Let $n$, $t_m$, $y_m$ and $\by_m$
be defined by the continued fraction expansion of $p/k$.
Then the Takahashi lengths $l_{j+1}$ and truncated
Takahashi lengths $\bl_{j+1}$
are defined as
\begin{equation}\label{taka}
\left.
\begin{aligned}
l_{j+1}&=  y_{m-1}+(j-t_m)  y_m \\
\bl_{j+1}&=\by_{m-1}+(j-t_m)\by_m
\end{aligned} \right\}
\; \text{ for } \; 0\leq m \leq n,~t_m<j\leq t_{m+1}+\delta_{n,m}.
\end{equation}
\end{definition}

Finally, vectors $\vQ^{(j)}$ $(j=1,\ldots,t_{n+1}+1)$ are needed to
specify parities of summation variables.
For $1\leq i\leq t_{n+1}$ and $0\leq m\leq n$ such that 
$t_m<j\leq t_{m+1}+\delta_{m,n}$ 
the components of $\vQ^{(j)}$ are recursively defined as
\begin{equation}\label{eqn_Qij}
Q^{(j)}_{i}=\begin{cases} 0 & \text{for $j \leq i \leq t_{n+1}$},\\  
j-i & \text{for $t_m\leq i<j$},\\
Q^{(j)}_{i+1}+Q^{(j)}_{t_{m'}+1} 
& \text{for $t_{m'-1}\leq i<t_{m'},~1\leq m'\leq m$},
\end{cases}
\end{equation}
with $Q_{t_n+1}^{(t_n+1)}=0$ for $\nu_n=0$.

In the next section we repeatedly use the above definitions
and it will be convenient to refer to all of the equations of this section
as ``the Takahashi--Suzuki (TS) decomposition of $p/k$''~\cite{Takahashi}.

\subsection{The supernomial boson--fermion identities}\label{sec_FB}
We are now prepared for the polynomial boson--fermion identities 
based on the continued fraction expansion of $p/k$, and 
depending on multiple finitization parameters.
To this end let us first define 
the generalizations of the bosonic and fermionic side 
of~\eqref{dualAG}, respectively. The bosonic side involves
the supernomials $T(\vL,a)$ of equation~\eqref{duality}.

\begin{definition}\label{def_B}
Consider the TS decomposition of $p/k$ with $p$ and $k$ positive
integers and $\gcd(p,k)$=1.
Fix a positive integer $N$ such that $N<(p-1)/k-1$ and 
let $\vL\in \Integer^N$, $\vl=T^{-1}\vL$.
Choose integers $a$ and $b$ such that $a+b+\l_N$ is even and such that 
$a=l_{\alpha+1}$ and $b=l_{\beta+1}\geq 2$ are Takahashi lengths. Then 
\begin{align}\label{eqn_B}
B_{a,b}^{(p,k,N)}(\vL)=q^{\frac{1}{4N}(b-a)^2}
\sum_{j=-\infty}^{\infty}\Bigl\{
q^{\frac{j}{N}\left(p(p-kN)j+p r-(p-kN)a\right)}
T\bigl(\vL,\tfrac{b-a}{2}+pj\bigr) &\notag\\
-q^{\frac{1}{N}(pj+a)((p-kN)j+r)}
T\bigl(\vL,\tfrac{b+a}{2}+pj\bigr)&\Bigr\},
\end{align}
with $r=b-N(b-\bar{b})$ and $\bar{b}=\bl_{\beta+1}$.
\end{definition}

\begin{definition}\label{def_F}
Fix all parameters as in definition~\ref{def_B}. Then
\begin{equation}\label{eqn_F}
F_{a,b}^{(p,k,N)}(\vL)=q^{\Delta_{a,b}}
\sum_{\substack{\vm\in \Integer_+^{\,t_{n+1}} \\ 
\vm\equiv \vQ_{a,b} \pmod{2}}}
q^{\frac{1}{4}\vm B\vm-\frac{1}{2}\vA_{a,b}\vm}
\prod_{j=1}^{t_{n+1}} 
\qbin{m_j+n_j}{m_j}.
\end{equation}
Here the $(\vm,\vn)$-system is given by
\begin{equation}\label{eqn_mn}
\vm+\vn=\frac{1}{2}\Bigl(\I_B\vm+\vu_a+\vu_b+\sum_{i=1}^N L_i\ve_i\Bigr),
\end{equation}
where 
\begin{align}\label{eqn_u}
\vu_a&=\ve_{\alpha}-\sum_{i=m+1}^n \ve_{t_i}\quad \text{for }
t_m <\alpha \leq t_{m+1}+\delta_{m,n}, \notag \\
\vu_b&=\ve_{\beta}-\sum_{i=m'+1}^n \ve_{t_i} \quad \text{for }
t_{m'} <\beta \leq t_{m'+1}+\delta_{m',n}.
\end{align}
In addition, for $t_i<j\leq t_{i+1}$
$(i=0,\ldots,n)$,
\begin{equation}
(\vA_{a,b})_j=\begin{cases} (\vu_b)_j & \text{for $i$ odd},\\
(\vu_a)_j & \text{for $i$ even}. \end{cases}
\end{equation}
The restriction $\vm\equiv \vQ_{a,b} \pmod{2}$ on the sum is an abbreviation 
of $m_i\equiv \left(\vQ_{a,b}\right)_i \pmod{2}$ for all 
$i=1,\ldots,t_{n+1}$. Here
\begin{equation}\label{eqn_Q}
\vQ_{a,b}=\sum_{j=1}^{t_{n+1}} (\vu_a+\vu_b)_j\vQ^{(j)}
+(\delta_{\alpha,t_{n+1}+1}+\delta_{\beta,t_{n+1}+1})\vQ^{(t_{n+1}+1)}
+\sum_{j=1}^N L_j\vQ^{(j)},
\end{equation}
with vectors $\vQ^{(j)}$ defined in equation~\eqref{eqn_Qij}.
Finally, $\Delta_{a,b}$ is fixed through the condition
\begin{equation}
\left.q^{-\frac{1}{4N}(b-a)^2}F_{a,b}^{(p,k,N)}(\vL)\right|_{q=0}=1
\qquad \text{for $\ell_N\geq |b-a|$ and $\vL\in\Integer^N_+$}.
\end{equation}
\end{definition}

Before proceeding, let us remark that the parity restriction 
$\vm \equiv \vQ_{a,b} \pmod{2}$ can be derived from 
$m_{t_{n+1}}\equiv \delta_{\alpha,t_{n+1}+1}+\delta_{\beta,t_{n+1}+1}
\pmod{2}$ and the condition that $\vn$ in~\eqref{eqn_mn} is a vector
with integer components.

For $\vL=L\ve_1$ $(L\geq 0)$, the bosonic function in \eqref{eqn_B} 
corresponds (up to a prefactor) to the one-dimensional configuration sum of 
the Andrews--Baxter--Forrester (ABF) model~\cite{ABF84} as obtained in 
ref.~\cite{FB85}. This follows from equation \eqref{T_L1}, yielding  
\begin{align}\label{FB}
B_{a,b}^{(p,k,N)}(L\ve_1)
=q^{\frac{1}{4}(b-a)^2}\sum_{j=-\infty}^{\infty}\biggl\{
q^{j(p(p-k)j+p\bar{b}-(p-k)a)}\qbin{L}{\frac{L+b-a}{2}+pj}&\notag\\
-q^{(pj+a)((p-k)j+\bar{b})}\qbin{L}{\frac{L+b+a}{2}+pj}&\biggr\}.
\end{align}
The fermionic function in \eqref{eqn_F} for $\vL=L\ve_1$ corresponds 
to the Berkovich--McCoy finitization~\cite{BM96} of the Virasoro characters of 
the minimal model $M(p-k,p)$.

Actually, in ref.~\cite{FB85} two functions $D_L(a,b,b+1)$ and $D_L(a,b,b-1)$
have been considered, with~\eqref{FB} corresponding to $D_L(a,b,b-1)$.
They are however not independent, but related by the symmetry
$D_L(a,b,b+1)=D_L(p-a,p-b,p-b-1)$ and hence it is sufficient 
to restrict attention to generalizations of $D_L(a,b,b-1)$.
We should also mention that in ref.~\cite{FB85} $D_L(a,b,b-1)$ is defined
for all $1\leq a\leq p-1$ and $2\leq b\leq p-1$ and not just for
$a,b$ being Takahashi lengths. The reason for restricting to
Takahashi lengths is that in the more general case 
the fermionic functions become rather complicated~\cite{BMS96}.

The following theorem claims a polynomial boson--fermion type relation.
\begin{theorem}\label{theo_FB}
For $b=l_{\beta+1}$ with $\beta\geq N$ and $\vL\in\Integer^N_+$, the functions
$B_{a,b}^{(p,k,N)}(\vL)$ and $F_{a,b}^{(p,k,N)}(\vL)$ of
definitions \ref{def_B} and \ref{def_F} satisfy the identity
\begin{equation}\label{eqn_FB}
B_{a,b}^{(p,k,N)}(\vL)=F_{a,b}^{(p,k,N)}(\vL).
\end{equation}
\end{theorem}

\begin{proof}
We show that both sides of \eqref{eqn_FB} satisfy the recursion relation
\begin{equation}\label{eqn_Xrec}
X(\vL)=X(\vL-2\ve_i)+q^{\frac{1}{2}(L_i-1)}X(\vL+\ve_{i-1}-2\ve_i+\ve_{i+1}),
\end{equation}
for $1\leq i\leq N-1$.
Since in ref.~\cite{BMS96} identity~\eqref{eqn_FB} has been shown to hold 
for $\vL=L\ve_1$ with $L\in \Integer_+$ and since 
$q^{\frac{1}{4}\vL T^{-1}\vL-\frac{a^2}{N}}X(\vL;1/q)$ satisfies~\eqref{Xrec},
this establishes theorem~\ref{theo_FB} by lemma~\ref{lem_complete}.

The recurrences \eqref{eqn_Trec} for the supernomials immediately imply 
that $B_{a,b}^{(p,k,N)}(\vL)$ satisfies \eqref{eqn_Xrec}. In order to show 
that $F_{a,b}^{(p,k,N)}(\vL)$ satisfies \eqref{eqn_Xrec}, we apply the
$q$-binomial recurrence \eqref{bi2} to the term $j=i$ of the product 
in~\eqref{eqn_F}, 
\begin{equation*}
\prod_{j=1}^{t_{n+1}}\qbin{m_j+n_j}{m_j}=
\prod_{j=1}^{t_{n+1}}\qbin{m_j+n_j-\delta_{i,j}}{m_j}+
q^{n_i}\prod_{j=1}^{t_{n+1}}\qbin{m_j+n_j-\delta_{i,j}}{m_j
-\delta_{i,j}}.
\end{equation*}
The first term on the right-hand side directly yields
$F_{a,b}^{(p,k,N)}(\vL-2\ve_i)$. In the second term we change 
the summation variable $m_i\to m_i+1$. 
Since $\vn$ can be expressed in terms of $\vm$ as
$$
n_j=\frac{1}{2}\left(-(B\vm)_j+(\vu_a)_j+(\vu_b)_j+\theta(j\leq N)L_j\right),
$$
this variable change gives
\begin{equation*}
\frac{1}{4}\vm B\vm-\frac{1}{2}\vA_{a,b}\vm+n_i\to
\frac{1}{2}(L_i-\frac{1}{2}B_{i,i})+\frac{1}{2}(\vu_a+\vu_b-\vA_{a,b})_i
+\frac{1}{4}\vm B\vm-\frac{1}{2}\vA_{a,b}\vm.
\end{equation*}
Now observe that the condition $N<(p-1)/k-1$ implies
that $N\leq \nu_0$, so that $i<\nu_0$, with $\nu_0$ defined 
in~\eqref{c_fraction}. This means that $B_{i,i}=2$ and
$(\vu_a+\vu_b-\vA_{a,b})_i=(\vu_b)_i=0$ since 
$\beta\geq N$. Hence we see that after the variable change
the second term yields
$q^{\frac{1}{2}(L_i-1)}F_{a,b}^{(p,k,N)}(\vL+\ve_{i-1}
-2\ve_i+\ve_{i+1})$ and \eqref{eqn_Xrec} is proven for $X=F$.
\end{proof}

{}From the proof of theorem~\ref{theo_FB} follows
that $\Delta_{a,b}$ in \eqref{eqn_F} is independent of $N$, and can
therefore be determined by considering the case $N=1$.\footnote{
Equation \eqref{dualAG} may seem to contradict this.
Note however that in theorem~\ref{theo_FB} we can choose $p,k$ and $N$
independently. Equation~\eqref{dualAG} simply corresponds to the
special choice $p=2N+3$, introducing $N$ dependence in ``$p$-dependent''
quantities.}
The explicit expression for $\Delta_{a,b}$ is quite involved and
can be found in ref.~\cite{BMS96}.

Finally notice that for $X=B,F$ and $N\leq M$
\begin{equation}\label{X_M}
X_{a,b}^{(p,k,M)}\bigl( (L_1,\dots,L_N,0,\dots,0)\bigr)=
X_{a,b}^{(p,k,N)}\bigl( (L_1,\dots,L_N)\bigr).
\end{equation}
For $X=B$ this identity is quite remarkable since, due to the fact that
the superscript on the left--hand side involves $M$ and the one on the
right--hand side involves $N$, the quadratic exponents of $q$ in the definition
of the bosonic side~\eqref{eqn_B} are different. Relation~\eqref{X_M} for 
$X=B$ may, however, be deduced from property~\eqref{T_M} of the 
$q$-supernomials.
For $X=F$ relation \eqref{X_M} follows directly from the 
$(\vm,\vn)$-system \eqref{eqn_mn} and the fact that $\Delta_{a,b}$ is 
$N$ independent. 

\subsection{The case $k=1$}\label{sec55}
The identities of section~\ref{sec_FB} are very general, but have the
drawback of being rather implicit. This warrants elaborating the simple
but important case $k=1$.

As noted before, for $k=1$ the TS decomposition is trivial, yielding $n=0$,
$p=\nu_0+3$ and $t_1=p-3$, so that $\I_B=\I$ and $B=C$ with $\I$ and $C$
the incidence matrix and Cartan matrix of A$_{p-3}$, respectively.
The (truncated) Takahashi lengths are given by
$l_{j+1}=j+1$ and $\bl_{j+1}=j$ for $0\leq j\leq p-2$.
This implies the following simplifications for some of the quantities
of definition~\ref{def_F}; $\vA_{a,b}=\vu_a=\ve_{a-1}$, $\vu_b=\ve_{b-1}$
and $\vQ_{a,b}=\vQ^{(a-1)}+\vQ^{(b-1)}+\sum_{i=2}^N L_i\vQ^{(i)}$ where
$\vQ^{(j)}=\ve_{j-1}+\ve_{j-3}+\cdots$.
Finally, we find that $\Delta_{a,b}=(b-a)/4$.
We can thus conclude the following supernomial identity $(m_0=0$):
\begin{multline}\label{unitary}
\sum_{j=-\infty}^{\infty}\Bigl\{
q^{\frac{j}{N}\left(p(p-N)j+p r-(p-N)a\right)}
T\bigl(\vL,\tfrac{b-a}{2}+pj\bigr)
-q^{\frac{1}{N}(pj+a)((p-N)j+r)}
T\bigl(\vL,\tfrac{b+a}{2}+pj\bigr)\Bigr\} \\
=q^{\frac{1}{4N}(b-a)(a-b+N)} \hspace{-3mm}
\sum_{\substack{\vm\in \Integer_+^{\,p-3} \\
\vm \equiv \vQ_{a,b} \pmod{2}}} \hspace{-3mm}
q^{\frac{1}{4}\vm C\vm-\frac{1}{2}m_{a-1}}
\prod_{j=1}^{p-3} \qbin{m_j+n_j}{m_j},
\end{multline}
for $\vL\in\Integer^N_+$, $1\leq a \leq p-1$, $N+1\leq b\leq p-1$ and $r=b-N$.
The corresponding $(\vm,\vn)$-system is given by
$$
\vm+\vn=\frac{1}{2}\Bigl(\I\vm+\ve_{a-1}+\ve_{b-1}+\sum_{i=1}^N L_i\ve_i\Bigr).
$$

For $\vL=L\ve_1$ $(L\geq 0)$, equation \eqref{unitary}, which has been proven 
in refs.~\cite{Berkovich94,W96a,W96b}, corresponds to an identity
for the one-dimensional configuration sums of the 
Andrews--Baxter--Forrester model~\cite{ABF84}. 
More generally, for $\vL=L\ve_N$, \eqref{unitary}
is an identity for configuration sums of RSOS models of
Date et al.~\cite{DJKMO87,DJKMO88},
and has first been proven in refs.~\cite{S96a,S96}.

\subsection{$q$-Series limits of theorem~\ref{theo_FB}}\label{sec_limits}
\subsubsection{Limits related to generalizations of the A$_1^{(1)}$ branching 
functions}
\label{sec_bf}
In this section we consider the limit $L_N\to\infty$ of the 
polynomials in~\eqref{eqn_B} and~\eqref{eqn_F},
and show how this limit relates to the branching functions of the 
A$_1^{(1)}$ cosets
\begin{equation}\label{coset}
\frac{({\rm A}^{(1)}_1)_N \times ({\rm A}^{(1)}_1)_{N'}}
{({\rm A}^{(1)}_1)_{N+N'}}\qquad N\in\Integer,~N'\in\Rational.
\end{equation}
We denote the normalized branching functions of these cosets by
$\hat{\chi}^{(P,P';N)}_{r,s;\ell}$ where
$$N'=\frac{N P}{P'-P}-2 \quad \text{or}\quad N'=-2-\frac{N P'}{P'-P}\, ,$$
with the restrictions $P<P'$, $P'-P \equiv 0 \pmod{N}$ and 
gcd$(\tfrac{P'-P}{N},P')=1$. 

Before presenting the $q$-series limits of the supernomial identities,
we first define generalizations of the A$_1^{(1)}$ branching functions.
\begin{definition}\label{def_bf}
Let $N,P,P'$ be positive integers such that $P<P'$, $P'-P \equiv 0 \pmod{N}$ 
and gcd$(\tfrac{P'-P}{N},P')=1$. Fix $\sigma=0,1$, $\vL\in\Integer_+^{N-1}$
and choose $1\leq r<P$ and $1\leq s<P'$ such that
$r-s+N(C^{-1}\vL)_{N-1}+N\sigma$ is even. Then
\begin{multline}\label{super_bf}
\hat{\chi}_{r,s;\vL,\sigma}^{(P,P';N)}(q)
=q^{-\frac{\vL C^{-1}\vL}{2(N+2)}} \\
\times
\sum_{0\leq m \leq N/2} c_{2m}^{\vL,\sigma}(q) \biggl(
\sum_{\substack{j \in \Integer\\ m_{r-s}(j)\equiv \pm m \pmod{N}}} 
\hspace{-6mm}  q^{\frac{j}{N}(jPP'+P'r-Ps)}
- \hspace{-4mm} \sum_{\substack{j\in \Integer \\ 
m_{r+s}(j)\equiv \pm m \pmod{N} }} 
\hspace{-6mm} q^{\frac{1}{N}(jP'+s)(jP+r)} \biggr),
\end{multline}
where $m_a(j)=(a/2+P'j)$.
The sum over $m$ runs over integers if $r-s$ is even and half-integers
if $r-s$ is odd.
The string-like function $c_{2m}^{\vL,\sigma}$ is defined 
in~\eqref{super_string}.
\end{definition}

Setting $\vL=\ve_{\ell}$ and using~\eqref{string},
this reduces to the normalized branching 
functions of the cosets~\eqref{coset}, in the representation 
of ref.~\cite{ACT91} (for the unitary case $P'=P+N$, so that
$N'\in \Integer$, see also refs.~\cite{BNY88,KMQ88,R88}),
\begin{equation}\label{bf_norm}
\hat{\chi}_{r,s;\ell}^{(P,P';N)}(q)=\begin{cases}
\hat{\chi}_{r,s;\ve_{\ell},0}^{(P,P';N)}(q) & \text{for $0\leq \ell<N$},\\[2mm]
\hat{\chi}_{r,s;{\bs 0},1}^{(P,P';N)}(q) & \text{for $\ell=N$}.
\end{cases}
\end{equation}

The generalized branching functions in definition~\ref{def_bf} arise as the 
following limit of the bosonic polynomials of equation~\eqref{eqn_B}. 
\begin{lemma}\label{lem_B_lim}
Fix $\sigma=0,1$ and positive integers $N,M,p,k$ such that 
$1\leq N\leq M<(p-1)/k-1$ and gcd$(p,k)=1$, and define the
TS decomposition of $p/k$.
Let $a,b,r$ and $B_{a,b}^{(p,k,M)}$ be as in definition~\ref{def_B},
and $\vL\in\Integer_+^{N-1}$ such that $r-a+N(C^{-1}\vL)_{N-1}+N\sigma$ 
is even. Then for $|q|<1$,
\begin{equation}\label{B_lim}
\hat{\chi}_{r,a;\vL,\sigma}^{(p-kN,p;N)}(q)=
q^{-\frac{1}{4N}(b-a)^2-\frac{1}{4}\vL C^{-1}\vL}
\lim_{\substack{L_N\to\infty\\L_N\equiv \sigma \pmod{2}}} 
B_{a,b}^{(p,k,M)}\left( (\vL,L_N,0,\ldots,0)\right).
\end{equation}
\end{lemma}

\begin{proof}
Because of~\eqref{X_M} with $X=B$, it is sufficient to establish~\eqref{B_lim}
for $N=M$. From corollary~\ref{cor_super_string} and the 
symmetries~\eqref{super_string_sym}, equation~\eqref{B_lim} follows 
immediately.
\end{proof}

Thanks to theorem~\ref{theo_FB} and lemma~\ref{lem_B_lim} we obtain the 
following fermionic representation for the extended branching functions.
\begin{corollary}\label{cor_ferm_bf}
Define all quantities as in lemma~\ref{lem_B_lim}. Then
\begin{equation*}
\hat{\chi}_{r,a;\vL,\sigma}^{(p-kN,p;N)}(q)=
q^{-\frac{1}{4N}(b-a)^2-\frac{1}{4}\vL C^{-1}\vL+\Delta_{a,b}}
\sum_{\substack{\vm\in\Integer_+^{t_{n+1}}\\ \vm\equiv \vQ_{a,b}\pmod{2}}}
\frac{q^{\frac{1}{4}\vm B\vm-\frac{1}{2}\vA_{a,b}\vm}}{(q)_{m_N}}
\prod_{\substack{j=1\\ j\neq N}}^{t_{n+1}}\qbin{m_j+n_j}{m_j},
\end{equation*}
with $(\vm,\vn)$-system 
$\vm+\vn=\frac{1}{2}(\I_B\vm+\vu_{a}+\vu_{b}+\sum_{k=1}^{N-1}L_k\ve_k)$,
$\vu_{a}$ and $\vu_{b}$ as in \eqref{eqn_u} and 
$\vQ_{a,b}$ as in \eqref{eqn_Q} with $L_N$ therein replaced by $\sigma$. 
\end{corollary}
For all but the cases where the above corresponds to the A$_1^{(1)}$
branching functions, we believe this result to be new.
For the special case of equation~\eqref{bf_norm}, 
the above was found in part in ref.~\cite{FQ94} and in general in
refs.~\cite{BM96,BMS96} for $N=1$, in 
ref.~\cite{S96a} for $N\geq 2$ and $k=1$, 
in ref.~\cite{W96} for $N\geq 2$, $k=2$ and $p=2N+3$ and in
ref.~\cite{BMSW} for general $N\geq 2$.

\subsubsection{Further limits}\label{sec_further_lim}
We briefly discuss the more general limits $L_{k_i}\to\infty$ 
$(i=1,\ldots,h)$ of theorem~\ref{theo_FB}, where 
$1\leq k_1<k_2<\cdots <k_h\leq N$. As in definition~\ref{def_b}, 
$K=\{ k_1,\ldots,k_h\}$ and $\bar{K}=\{1,\ldots,N\}-K$. In addition 
we denote $K'=\{1,\ldots,t_{n+1}\}-K$.  Then we obtain
\begin{equation*}
\lim_{\substack{L_{k_1},\ldots,L_{k_h}\to\infty\\
L_{k_i}\equiv \sigma_{k_i} \pmod{2},~(1\leq i\leq h)}}
F_{a,b}^{(p,k,N)}(\vL)=q^{\Delta_{a,b}}
\sum_{\substack{\vm\in\Integer_+^{t_{n+1}}\\ \vm\equiv \vQ_{a,b} \pmod{2}}}
\frac{q^{\frac{1}{4}\vm B\vm-\frac{1}{2}\vA_{a,b}\vm}}
{(q)_{m_{k_1}}\cdots (q)_{m_{k_h}}}
\prod_{j\in K'} \qbin{m_j+n_j}{m_j}
\end{equation*}
where $\vm+\vn=\frac{1}{2}(\I_B\vm+\vu_a+\vu_b+\sum_{k\in\bar{K}}
L_k\ve_k)$, $\vu_a$ and $\vu_b$ as in~\eqref{eqn_u} and $\vQ_{a,b}$ given by
$$
\vQ_{a,b}=\sum_{j=1}^{t_{n+1}} (\vu_a+\vu_b)_j\vQ^{(j)}
+(\delta_{\alpha,t_{n+1}+1}+\delta_{\beta,t_{n+1}+1})\vQ^{(t_{n+1}+1)}
+\sum_{j\in K} \sigma_j\vQ^{(j)}+\sum_{j\in \bar{K}} L_j\vQ^{(j)},
$$
recalling that $a=l_{\alpha+1}$ and $b=l_{\beta+1}$.

According to equation \eqref{limit2}, 
the same limit, but now for $B_{a,b}^{(p,k,N)}$, is obtained by 
replacing the $q$-supernomials in~\eqref{eqn_B} by 
$b_{2m}^{\{L_k|k\in\bar{K}\}\{\sigma_k|k\in K\}}$ with $m=(b-a)/2+pj$ for 
the first $q$-supernomial and $m=(b+a)/2+pj$ for the second $q$-supernomial.

A yet different set of $q$-series identities can be obtained from
theorem~\ref{theo_FB} by first transforming $q\to 1/q$ before 
taking some (or all) of the components of $\vL$ to infinity.
Recalling~\eqref{duality} we find that
\begin{equation*}
q^{\frac{1}{4}\vL T^{-1}\vL} B_{a,b}^{(p,k,N)}(\vL;1/q) 
=\sum_{j=-\infty}^{\infty}\biggl\{ 
q^{j\left(pkj+p(b-\bar{b})-ka\right)} \qbin{\vL}{\tfrac{b-a}{2}+pj}
-q^{(pj+a)(kj+b-\bar{b})} \qbin{\vL}{\tfrac{b+a}{2}+pj}\biggr\}, 
\end{equation*}
where we have used that $r=b-N(b-\bar{b})$.
Hence, by lemma~\ref{lem_limit1}, we see that for all $m=1,\dots,N$,
\begin{equation*}
\lim_{L_m \to \infty}
q^{\frac{1}{4}\vL T^{-1}\vL} B_{a,b}^{(p,k,N)}(\vL;1/q) 
= \frac{1}{(q)_{\infty}}
\sum_{j=-\infty}^{\infty}\Bigl\{
q^{j\left(pkj+p(b-\bar{b})-ka\right)}
-q^{(pj+a)(kj+b-\bar{b})} \Bigr\},
\end{equation*}
independent of $N$.
The right-hand side of this expression is recognized as the
bosonic form of the (normalized) Virasoro characters
$\hat{\chi}_{b-\bar{b},a}^{(k,p)}(q)$.

Of course we can carry out a similar calculation for the function $F$ of
theorem~\ref{theo_FB}, but since the resulting fermionic representations
for the characters $\hat{\chi}_{b-\bar{b},a}^{(k,p)}(q)$ coincide with those
given in refs.~\cite{BM96,BMS96}, we refer the interested reader to the
literature.

\section{Discussion}
We conclude with some final remarks about the results of this paper.

First we comment on the restriction $N<(p-1)/k-1$ imposed
on the $q$-supernomial identity of theorem~\ref{theo_FB}.
Recall the continued fraction expansion of $p/k$ as given in
equation~\eqref{c_fraction}. This expansion defines integers
$\nu_0,\dots,\nu_n$ which are used to define an $(\vm,\vn)$-system
of dimension $\nu_0+\cdots+\nu_n$.
As can be seen from equation \eqref{I_B} ,
the components of the $(\vm,\vn)$-system
with labels between $\nu_0+\cdots+\nu_i$ and $\nu_0+\cdots+\nu_{i+1}$ 
can be viewed as one unit. Each of these different
units is referred to as a ``zone''~\cite{Takahashi}.
Now the restriction in our theorem implies that $N \leq \nu_0$
and hence that our recurrences essentially act within the zeroth
zone of the $(\vm,\vn)$-system only. To obtain $q$-supernomial
identities which are true for $N>\nu_0$ one needs to generalize
the $q$-supernomials to satisfy the recurrences implied by the full
set of zones of the $(\vm,\vn)$-system.
Such modified $q$-supernomials would  thus reflect the
complete continued fraction expansion of $p/k$ and not
just its integer part.
Indeed polynomial boson--fermion identities (related to
the Gordon--G\"ollnitz identities) with finitization parameter $L$ 
beyond zone zero have been found in refs.~\cite{BMO96,BM95}.
It would be an interesting problem to generalize these to identities
depending on multiple finitization parameters, using more general
$q$-supernomials.

Another point of interest is the combinatorial interpretation
of the bosonic and fermionic polynomials of equations \eqref{eqn_B}
and \eqref{eqn_F}. 
For $\vL=L \ve_i$ ($i=1,\dots,N)$ these polynomials correspond to 
one-dimensional configuration sums of solvable lattice 
models~\cite{ABF84,FB85,DJKMO87,DJKMO88}.
In those cases we have an interpretation in terms of 
weighted lattice paths. Another interpretation, in terms of
partitions with prescribed hook differences, exists for 
$\vL=L\ve_1$~\cite{ABBBFV87}. 
For the general $N$-dimensional case, however, we have not been
able to generalize either of these.
What is possible is to give a simple interpretation
in terms of lattice paths when $q=1$. In fact, instead of
working with lattice paths, we give an equivalent matrix formulation
in the following.

Fix an integer $p\geq 4$ and define a family of incidence matrices
$A_0,\ldots,A_{p-2}$ as follows:
\begin{equation*}
A_k A_1 = A_{k-1}+A_{k+1} \quad \text{ for } k=1,\dots,p-3,
\end{equation*}
with $A_0=I$ the $(p-1)\times (p-1)$ identity matrix and $A_1=\I$
the incidence matrix of the Lie algebra A$_{p-1}$, 
with entries $\I_{a,b}=\delta_{|a-b|,1}$.
The above set of matrices forms a commuting family.
Some further notable properties are $A_{p-2}=Y$ with $Y_{a,b}=\delta_{a,p-b}$ 
and $A_k Y = A_{p-k-2}$. 
We now claim for $\vL\in \Integer_+^{p-3}$,
$a,b=1,\dots p-1$, such that $a+b+L_1+L_3+\cdots$ is even, 
that
\begin{equation*}
(A_1^{L_1}A_2^{L_2} \dots A_{p-3}^{L_{p-3}})_{a,b}
=\sum_{j=-\infty}^{\infty}\biggl\{
\binom{\vL}{\tfrac{b-a}{2}+pj} - 
\binom{\vL}{\tfrac{b+a}{2}+pj}\biggr\},
\end{equation*}
with on the right-hand side the 
supernomials of definition~\ref{def_supernomials} (with $N=p-3$).

Finally we note that the
boson--fermion identities for the branching functions of the
A$_1^{(1)}$ coset theories~\eqref{coset}, which occur as a particular limit
of the polynomial identities of theorem~\ref{theo_FB},
can also be derived
using the higher-level Bailey lemma of refs.~\cite{SW96a,SW96b}, as
is the subject of ref.~\cite{BMSW}.
It is intriguing to note that though
the method of ref.~\cite{BMSW} and that of the present paper
are very different, both
rely on the use of the polynomial identities for the minimal models
$M(p-k,p)$~\cite{BM96,BMS96}.
Whereas in ref.~\cite{BMSW} they provide the necessary Bailey pairs
as input for the higher-level lemma, in this paper
they serve as initial conditions for the recursive proofs
of the polynomial identities.

\section*{Acknowledgements}
We are thankful to Barry McCoy for discussions and 
to Alexander Berkovich for drawing our attention to 
reference~\cite{Bressoud92}. This work was 
supported by the Australian Research Council and NSF grant DMR9404747.

\section*{Note added}
We now have obtained a combinatorial interpretation of the boson-fermion
identities of theorem~\ref{theo_FB} in terms of the inhomogeneous lattice paths
introduced by Nakayashiki and Yamada~\cite{NY95} and
Lascoux, Leclerc and Thibon~\cite{LLT95}.


\begin{thebibliography}{99}

\bibitem{ACT91}
C.~Ahn, S.-W.~Chung and S.-H.~Tye,
{\em New parafermion, SU(2) coset and $N=2$ superconformal field theories},
Nucl.\ Phys.\ B {\bf 365} (1991), 191--240.

\bibitem{Andrews70}
G.~E.~Andrews,
{\em A polynomial identity which implies the Rogers--Ramanujan identities},
Scripta Math.\ {\bf 28} (1970), 297--305.

\bibitem{Andrews72}
G.~E.~Andrews,
{\em Sieves in the theory of partitions},
Amer.\ J. Math.\ {\bf 94} (1972), 1214--1230.

\bibitem{Andrews74}
G.~E.~Andrews,
{\em An analytic generalization of the Rogers--Ramanujan identities
for odd moduli},
Prod.\ Nat.\ Acad.\ Sci.\ USA {\bf 71} (1974), 4082--4085.

\bibitem{Andrews76}
G.~E.~Andrews,
{\em The Theory of Partitions},
Encyclopedia of Mathematics, vol.~2
(Addison-Wesley, Reading, Massachusetts, 1976).

\bibitem{Andrews79}
G.~E.~Andrews,
{\em Partitions and Durfee dissection},
Amer.\ J. Math.\ {\bf 101} (1979), 735--742.

\bibitem{Andrews85}
G.~E.~Andrews,
{\em $q$-Series: Their development and application in analysis, number
theory, combinatorics, physics, and computer algebra}, 
in CBMS Regional Conf.\ Ser.\ in Math.\ {\bf 66} 
(AMS, Providence, Rhode Island, 1985).

\bibitem{Andrews94}
G.~E.~Andrews,
{\em Schur's theorem, Capparelli's conjecture and $q$-trinomial
coefficients},
Contemp.\ Math.\ {\bf 166} (1994), 141--154.

\bibitem{AB87}
G.~E.~Andrews and R.~J.~Baxter,
{\em Lattice gas generalization of the hard
hexagon model. III. $q$-trinomial coefficients},
J. Stat.\ Phys.\ {\bf 47} (1987), 297--330.

\bibitem{ABBBFV87}
G.~E.~Andrews, R.~J.~Baxter, D.~M.~Bressoud, W.~H.~Burge, P.~J.~Forrester
and G.~Viennot,
{\em Partitions with prescribed hook differences},
Europ.\ J. Combinatorics {\bf 8} (1987), 341--350.

\bibitem{ABF84}
G.~E.~Andrews, R.~J.~Baxter and P.~J.~Forrester,
{\em Eight-vertex SOS model and generalized 
Rogers--Ramanujan--type identities},
J. Stat.\ Phys.\ {\bf 35} (1984), 193--266.

\bibitem{BNY88}
J.~Bagger, D.~Nemeschansky and S.~Yankielowicz,
{\em Virasoro algebras with central charge $c>1$},
Phys.\ Rev.\ Lett.\ {\bf 60} (1988), 389--392.

\bibitem{Berkovich94}
A.~Berkovich,
{\em Fermionic counting of RSOS states and Virasoro character formulas 
for the unitary minimal series $M(\nu,\nu+1)$: Exact results},
Nucl.\ Phys. B {\bf 431} (1994), 315--348.

\bibitem{BM96}
A.~Berkovich and B.~M.~McCoy,
{\em Continued fractions and fermionic representations for characters
of $M(p,p')$ minimal models},
Lett.\ Math.\ Phys.\ {\bf 37} (1996), 49--66.

\bibitem{BM95}
A.~Berkovich and B.~M.~McCoy,
{\em Generalizations of the Andrews--Bressoud identities for the $N=1$
superconformal model $SM(2,4\nu)$},
preprint BONN-TH-95-15, ITP-SB-95-29, hep-th/9508110. To appear in
Int.\ J. of Math.\ and Comp.\ Modelling.

\bibitem{BMO96}
A.~Berkovich, B.~M.~McCoy and W.~P.~Orrick,
{\em Polynomial identities, indices, and duality for the $N=1$ 
superconformal model $SM(2,4\nu)$},
J. Stat.\ Phys.\ {\bf 83} (1996), 795--837.

\bibitem{BMS96}
A.~Berkovich, B.~M.~McCoy and A.~Schilling,
{\em Rogers--Schur--Ramanujan type identities for the $M(p,p')$ minimal models
of conformal field theory},
preprint BONN-TH-96-07, ITP-SB-96-35, q-alg/9607020. 
To appear in Commun.\ Math.\ Phys.

\bibitem{BMSW}
A.~Berkovich, B.~M.~McCoy, A.~Schilling and S.~O.~Warnaar,
{\em Bailey flows and Bose--Fermi identities for the conformal
coset models  $({\rm A}^{(1)}_1)_N\times ({\rm A}^{(1)}_1)_{N'}/
({\rm A}^{(1)}_1)_{N+N'}$}, 
Nucl. Phys. B {\bf 499} [PM] (1997), 621--649.

\bibitem{Bressoud89}
D.~M.~Bressoud,
{\em In the land of OZ},
in $q$-Series and partitions,
ed., D.~Stanton, IMA Volume in Mathematics and its Applications 
(Springer, New York, 1989), 45--66.

\bibitem{Bressoud92}
D.~M.~Bressoud,
{\em Unimodality of Gaussian polynomials},
Discrete Math.\ {\bf 99} (1992), 17--24.

\bibitem{DJKMO87}
E.~Date, M.~Jimbo, A.~Kuniba, T.~Miwa and M.~Okado,
{\em Exactly solvable SOS models: Local height probabilities and 
theta function identities},
Nucl.\ Phys.\ B {\bf 290} [FS20] (1987), 231--273.

\bibitem{DJKMO88}
E.~Date, M.~Jimbo, A.~Kuniba, T.~Miwa and M.~Okado,
{\em Exactly solvable SOS sodels II: Proof of the star-triangle relation
and combinatorial identities},
Adv.\ Stud.\ in Pure Math.\ {\bf 16} (1988), 17--122.

\bibitem{FF82}
B.~L.~Feigin and D.~B.~Fuchs,
{\em Verma modules over the Virasoro algebra},
Topology (Leningrad, 1982), 230--245,
Lecture Notes in Math.\ {\bf 1060} (Springer, Berlin--New York, 1984).

\bibitem{FF83}
B.~L.~Feigin and D.~B.~Fuchs,
{\em Verma modules over the Virasoro algebra},
Funct.\ Anal.\ Appl.\ {\bf 17} (1983), 241--242.

\bibitem{FQ95}
O.~Foda and Y.-H.~Quano,
{\em Polynomial identities of the Rogers--Ramanujan type},
Int.\ J. Mod.\ Phys.\ A {\bf 10} (1995), 2291--2315.

\bibitem{FQ94}
Foda and Y.-H.~Quano,
{\em Virasoro character identities from the Andrews--Bailey construction},
Int.\ J. Mod.\ Phys.\ A {\bf 12} (1997), 1651--1676.

\bibitem{FB85}
P.~J.~Forrester and R.~J.~Baxter,
{\em Further exact solutions of the eight-vertex SOS model and 
generalizations of the Rogers--Ramanujan identities},
J. Stat.\ Phys.\ {\bf 38} (1985), 435--472.

\bibitem{GR90}
G.~Gasper and M.~Rahman,
{\em Basic Hypergeometric Series},
Encyclopedia of Mathematics, vol. 35 (Cambridge University Press, 
Cambridge, 1990).

\bibitem{JM84}
M.~Jimbo and T.~Miwa,
{\em Irreducible decomposition of fundamental modules for A$_l^{(1)}$
and C$_l^{(1)}$, and Hecke modular forms},
Adv.\ Stud.\ in Pure Math.\ {\bf 4} (1984), 97--119.  

\bibitem{KP84}
V.~G.~Kac and D.~H.~Peterson,
{\em Infinite-dimensional Lie algebras, theta functions and modular forms},
Adv.\ in Math.\ {\bf 53} (1984), 125--264.

\bibitem{KMQ88}
D.~Kastor, E.~Martinec and Z.~Qiu,
{\em Current algebra and conformal series},
Phys.\ Lett.\ B {\bf 200} (1988), 434--440.

\bibitem{KKMM93a}
R.~Kedem, T.~R.~Klassen, B.~M.~McCoy and E.~Melzer,
{\em Fermionic quasi-particle representations for characters of
$(G^{(1)})_1\times (G^{(1)})_1/(G^{(1)})_2$},
Phys.\ Lett.\ B {\bf 304} (1993), 263--270.

\bibitem{KKMM93b}
R.~Kedem, T.~R.~Klassen, B.~M.~McCoy and E.~Melzer,
{\em Fermionic sum representations for conformal field theory characters},
Phys.\ Lett.\ B {\bf 307} (1993), 68--76.

\bibitem{Kirillov95}
A.~N.~Kirillov,
{\em Dilogarithm identities},
Prog.\ Theor.\ Phys.\ Suppl.\ {\bf 118} (1995), 61--142.

\bibitem{LLT95}
A.~Lascoux, B.~Leclerc and J.-Y.~Thibon,
{\em Crystal graphs and $q$-analogues of weight multiplicities for the root
system A$_n$},
Lett. Math. Phys. {\bf 35} (1995), 359--374.

\bibitem{MacMahon}
P.~A.~MacMahon,
{\em Combinatory Analysis} vol. 2
(Cambridge University Press, Cambrigde, 1916).

\bibitem{NY95}
A.~Nakayashiki and Y.~Yamada,
{\em Kostka polynomials and energy functions in solvable lattice models},
q-alg/9512027. To appear in Sellecta Mathematica.

\bibitem{R88}
F.~Ravanini,
{\em An infinite class of new conformal field theories with extended algebras},
Mod.\ Phys.\ Lett.\ A {\bf 3} (1988), 397--412.

\bibitem{S96a}
A.~Schilling,
{\em Polynomial fermionic forms for the branching functions
of the rational coset conformal field theories 
$\widehat{su}(2)_M \times \widehat{su}(2)_N / \widehat{su}(2)_{M+N}$},
Nucl.\ Phys.\ B {\bf 459} (1996), 393--436.

\bibitem{S96}
A.~Schilling,
{\em Multinomials and polynomial bosonic forms for the branching functions
of the $\widehat{su}(2)_M \times \widehat{su}(2)_N / \widehat{su}(2)_{N+M}$
conformal coset models},
Nucl.\ Phys.\ B {\bf 467} (1996), 247--271.

\bibitem{SW96a}
A.~Schilling and S.~O.~Warnaar,
{\em A higher-level Bailey lemma},
Int.\ J. Mod.\ Phys.\ B {\bf 11} (1997), 189--195.

\bibitem{SW96b}
A.~Schilling and S.~O.~Warnaar,
{\em A higher-level Bailey lemma: Proof and application},
preprint ITP-SB-96-33, University of Melbourne No. 06-96, q-alg/9607014.
To appear in the Ramanujan Journal.

\bibitem{Schur17}
I.~J.~Schur,
{\em Ein Beitrag zur additiven Zahlentheorie und zur Theorie 
der Kettenbr\"uche},
S.-B.\ Preuss.\ Akad.\ Wiss.\ Phys.-Math.\ Kl.\ (1917), 302--321.

\bibitem{Takahashi}
M.~Takahashi and M.~Suzuki,
{\em One-dimensional anisotropic Heisenberg model at finite temperatures},
Prog.\ Theor.\ Phys.\ {\bf 48} (1972), 2187--2209.

\bibitem{W96a}
S.~O.~Warnaar,
{\em Fermionic solution of the Andrews--Baxter--Forrester model. 
I. Unification of TBA and CTM methods},
J. Stat.\ Phys.\ {\bf 82} (1996), 657--685.

\bibitem{W96b}
S.~O.~Warnaar,
{\em Fermionic solution of the Andrews--Baxter--Forrester model. 
II. Proof of Melzer's polynomial identities},
J. Stat.\ Phys.\ {\bf 84} (1996), 49--83.
    
\bibitem{W96}
S.~O.~Warnaar,
{\em The Andrews--Gordon identities and $q$-multinomial coefficients},
Commun.\ Math.\ Phys. {\bf 184} (1997), 203--232.

\end{thebibliography}
\end{document}